\documentclass[12pt,a4paper,reqno]{amsart}
\usepackage[ansinew]{inputenc}
\usepackage{amsthm,amsmath,amssymb}
\usepackage{enumerate}
\RequirePackage[colorlinks,citecolor=blue,urlcolor=blue]{hyperref}

\newcommand\NoBlackBoxes{\global\overfullrule0pt}
\NoBlackBoxes
\theoremstyle{plain} % default
\def\4{\kern1pt}

\def\6{\vphantom0}

\def\8{\kern-10pt}
\def\7#1{_{(#1)}}

\makeatletter
\let\serieslogo@\relax
\let\@setcopyright\relax

\def\speciallabelmark#1{\def\@currentlabel{#1}}
\makeatother

\begin{document}

\def\ffrac#1#2{\raise.5pt\hbox{\small$\4\displaystyle\frac{\,#1\,}{\,#2\,}\4$}}
\def\ovln#1{\,{\overline{\!#1}}}
\def\ve{\varepsilon}
\def\kar{\beta_r}
\def\theequation{\thesection.\arabic{equation}}
\def\E{{\mathbb E}}
\def\R{{\mathbb R}}
\def\C{{\mathbb C}}
\def\Z{{\mathbb Z}}
\def\P{{\mathbb P}}
\def\H{{\rm H}}
\def\Im{{\rm Im}}
\def\Tr{{\rm Tr}}

\def\k{{\kappa}}
\def\M{{\cal M}}
\def\Var{{\rm Var}}
\def\Ent{{\rm Ent}}
\def\O{{\rm Osc}_\mu}

\def\ep{\varepsilon}
\def\phi{\varphi}
\def\F{{\cal F}}
\def\L{{\cal L}}

\def\s{{\mathfrak s}}

\def\be{\begin{equation}}
\def\en{\end{equation}}
\def\bee{\begin{eqnarray*}}
\def\ene{\end{eqnarray*}}

\title[Esscher Transform and CLT]{ESSCHER TRANSFORM AND \\
THE CENTRAL LIMIT THEOREM
}

\author{Sergey G. Bobkov$^{1}$}
\thanks{1) 
School of Mathematics, University of Minnesota, Minneapolis, MN, USA
}

\author{Friedrich G\"otze$^{2}$}
\thanks{2) Faculty of Mathematics,
Bielefeld University, Germany
}

\begin{abstract}
The paper is devoted to the investigation of Esscher's transform on high dimensional 
Euclidean spaces in the light of its application to the central limit theorem. 
With this tool, we explore necessary and sufficient conditions of normal approximation
for normalized sums of i.i.d. random vectors in terms of the R\'enyi divergence of 
infinite order, extending recent one dimensional results.
\end{abstract}


\maketitle

\section{{\bf Introduction}}
\setcounter{equation}{0}

\vskip2mm
\noindent
Introduce the normalized sums
$$
Z_n = \frac{X_1 + \dots + X_n}{\sqrt{n}}
$$
of independent copies of a random vector $X$ in $\R^d$ with mean zero
and identity covariance matrix $I_d$. It is well known that, if $Z_n$ have densities
$p_n$ for large $n$, their distributions are convergent in total variation norm
to the standard normal law on $\R^d$ with density 
$\varphi(x) = \frac{1}{(2\pi)^{d/2}}\,\exp(-|x|^2/2)$. That is, we have
the convergence in $L^1$-norm for densities
$$
\int_{\R^d} |p_n - \varphi|\,dx \rightarrow 0 \quad
(n \rightarrow \infty),
$$
which was first emphasized by Prokhorov \cite{Pr}.
A much stronger property, which may or may not hold in general, 
is described by means of the R\'enyi divergence
$$
D_\alpha(p_n||\varphi) = \frac{1}{\alpha - 1} \log 
\int_{\R^d} \Big(\frac{p}{\varphi}\Big)^\alpha\,\varphi\, dx
$$
of order $\alpha \geq 1$ (the relative $\alpha$-entropy).
These distance-like functionals are increasing for growing $\alpha$, and 
in the limit they define the R\'enyi divergence of infinite order
$$
D_\infty(p_n||\varphi)  = {\rm ess\,\sup}_x\, \log\, \frac{p_n(x)}{\varphi(x)}.
$$

The convergence in relative entropy $D=D_1$ (which is the Kullback-Leibler distance)
was the subject of numerous investigations starting with Linnik \cite{L},
who initiated an information-theoretic approach to the central limit theorem. 
Let us only mention the work by Barron \cite{Bar} (for necessary and 
sufficient conditions) and later \cite{A-B-B-N}, \cite{B-C-G1}, \cite{B-C-G2} 
(for the problems of rates and Berry-Esseen bounds). The case $\alpha>1$ and,
in particular, the Pearson $\chi^2$-distance was treated in detail in \cite{B-C-G3}, 
where the following statement was obtained as a consequence of 
a certain characterization 
of the convergence with respect to $D_\alpha$ in terms of the Laplace transform.

\vskip5mm
{\bf Theorem 1.1.} {\sl Assume that $D_\alpha(p_n||\varphi) <\infty$
for every $\alpha \in (1,\infty)$ with some $n = n_\alpha$. For
the convergence $D_\alpha(p_n||\varphi) \rightarrow 0$ for any finite $\alpha$,  
it is necessary and sufficient that 
\be
\E\, e^{\left<t,X\right>} \leq e^{|t|^2/2} \quad {\sl for \ all} \, \, t \in \R^d.
\en
}

The assumption of finiteness of $D_\alpha(p_n||\varphi)$
may be equivalently stated as the property that $Z_n$ have bounded
densities for large $n$.

The inequality (1.1) describes a remarkable class of probability distributions
which appear naturally in many mathematical problems. In modern literature, 
(1.1) is often called strict sub-Gaussianity. We refer an interested reader 
to \cite{B-C-G4} for the history, references, and recent developments 
towards the problem of characterization of such distributions in dimension one.

The convergence in $D_\infty$ is equivalent to the convergence with respect to
$$
T_\infty(p_n||\varphi) = {\rm ess\,\sup}_x\, \frac{p_n(x) - \varphi(x)}{\varphi(x)}
$$
in view of the relation $T_\infty = e^{D_\infty} - 1$. This quantity looks
more natural than $D_\infty$ for local limit theorems, where 
$|p_n(x) - \varphi(x)|$ is commonly estimated on $\R^d$ with polynomially
growing weights. However, if $X$ is not normal, the statement
\be
{\rm ess\,\sup}_x\,  \frac{p_n(x) - \varphi(x)}{\varphi(x)} \rightarrow 0 \quad
{\rm as} \ n \rightarrow \infty,
\en
cannot be strengthened to absolute value under the supremum. Indeed, if
for example, $X$ is bounded, then $p_n(x)$ is compactly supported, and
the above ratio is equal to $-1$ for large $|x|$
(see also Corollary 6.5 and the relation (6.5) after it).

One of the purposes of this paper is to give necessary and sufficient conditions
for the multidimensional CLT such as (1.2) in terms of the Laplace transform
$L(t) = \E\,e^{\left<t,X\right>}$. Introduce $K(t) = \log L(t)$ (which 
is a convex smooth function) and define the function
$$
A(t) = \frac{1}{2}\, |t|^2 - K(t), \quad t \in \R^d.
$$
As before,
suppose that $(X_k)_{k \geq 1}$ are independent copies of the random
vector $X$ in $\R^d$ with mean zero and identity covariance matrix. Below
we assume that:

\vskip3mm
1) $Z_n$ has density $p_n$ with $T_\infty(p_n||\varphi) < \infty$ for some $n=n_0$;

2) $X$ is strictly sub-Gaussian, that is, $A(t) \geq 0$ for all $t \in \R^d$.

\vskip5mm
{\bf Theorem 1.2.} {\sl For the convergence $T_\infty(p_n||\varphi) \rightarrow 0$,
it is necessary and sufficient that the following two conditions are fulfilled:

\vskip3mm
$a)$ \ $A''(t) = 0$ for every point $t \in \R^d$ such that $A(t) = 0$;

$b)$ \ $\lim_{k \rightarrow \infty} A''(t_k) = 0$ for every sequence
$|t_k| \rightarrow \infty$ such that $A(t_k) \rightarrow 0$ as $k \rightarrow \infty$.
}

\vskip5mm
Here and elsewhere $A''$ denotes the Hessian, that is, the $d \times d$ 
matrix of second order partial derivatives of $A$.
The conditions $a)-b)$ may be combined in the requirement 
$$
\lim_{A(t) \rightarrow 0} A''(t) = 0, \ \ {\rm or \ equivalently} \ \
\lim_{A(t) \rightarrow 0} K''(t) = I_d
$$
which is a kind of continuity of $A''$ with respect to $A$. As we will see,
these conditions may also be stated in a formally weaker form as

\vskip4mm
$a')$ ${\rm det} K''(t) = 1$ for every $t \in \R^d$ such that $A(t) = 0$;

$b')$ $\lim_{k \rightarrow \infty} {\rm det} K''(t_k) = 1$ for any sequence
$|t_k| \rightarrow \infty$ such that $A(t_k) \rightarrow 0$ as $k \rightarrow \infty$.

\vskip4mm
In dimension $d=1$ Theorem 1.2 has been obtained in \cite{B-G1}, 
where $a)-b)$ are stated as a weaker condition
$\limsup_{A(t) \rightarrow 0} A''(t) \leq 0$.
The multidimensional situation turns out to be more complicated, since
it requires a careful treatment of eigenvalues of the matrix
$K''(t)$, when $A(t)$ approaches zero. Another ingredient 
in the proof is a quantitative version of the uniform local limit theorem,
which was recently developed~in \cite{B-G2}.

Assuming the strict sub-Gaussianity (1.1), the conditions $a)-b)$ may or may not
hold in general. This shows that the convergence in $D_\infty$ is stronger
than the convergence in $D_\alpha$ simultaneously for all finite $\alpha$.
Nevertheless, for a wide class of strictly sub-Gaussian distributions 
the Laplace transform possesses a separation-type property
\be
\sup_{|t| \geq t_0} \big[e^{-|t|^2/2}\, \E\,e^{\left<t,X\right>} \big] < 1 
\quad {\rm for \ all} \ t_0>0,
\en
This is a strengthened form of condition 2), which entails
properties $a)-b)$.

\vskip5mm
{\bf Corollary 1.3.} {\sl If a random vector $X$ with mean zero and 
identity covariance matrix satisfies $(1.3)$,
then $T_\infty(p_n||\varphi) \rightarrow 0$ as $n \rightarrow \infty$.
}

\vskip4mm
On the other hand, the case of equality in the sub-Gaussian bound (1.1) is
quite possible, and one can observe new features in the multidimensional 
case. While in dimension one, an equality $L(t) = e^{|t|^2/2}$ is only 
possible for a discrete set of points $t$, in higher dimensions the set of 
points where this equality holds may have dimension $d-1$. In order to 
clarify this behaviour, we will discuss the class of Laplace transforms which 
contain periodic components separately. Specializing Theorem 2.1 
to this class, the general characterization may be simplified.

\vskip5mm
{\bf Corollary 1.4.} {\sl Suppose that the function $\Psi(t) = L(t)\,e^{-|t|^2/2}$
is $h$-periodic for some vector $h \in \R^d_+$ $(h \neq 0)$. For the convergence
$T_\infty(p_n||\varphi) \rightarrow 0$ as $n \rightarrow \infty$,
it is necessary and sufficient that, for every $t \in [0,h]$,
$$
\Psi(t) = 1 \, \Rightarrow \, \Psi''(t) = 0.
$$
}

Our approach to Theorem 1.2 is based on the application of the Esscher
transforms $Q_h$ on $\R^d$ defined below. Recall that a random vector $X$ 
in $\R^d$ is called sub-Gaussian, if $\E\,e^{c|X|^2} < \infty$ for some $c>0$.

\vskip3mm
{\bf Definition 1.5.} Let $X$ be a sub-Gaussian random vector in $\R^d$ with 
distribution $\mu$. Introduce the family of probability distributions
$\mu_h = Q_h \mu$ on $\R^d$ with parameter $h \in \R^d$ which have
densities with respect to $\mu$
\be
\frac{d\mu_h(x)}{d\mu(x)} = \frac{1}{L(h)}\,e^{\left<h,x\right>}, \quad x \in \R^d.
\en

The early history of this transform in dimension one goes back to Esscher \cite{E} 
in actuarial science, Khinchin \cite{K} in statistical mechanics, and
Daniels \cite{D} who applied it to develop asymptotic expansions for densities. 

As a result, in (1.4) we obtain a semi-group of probability measures
$\{Q_h \mu\}_{h \in \R^d}$ on $\R^d$ with the following remarkable 
property: Every $Q_h$ transforms the convolution of several measures to the
convolutions of their $Q_h$-transforms. This is analogous to the property that
the Fourier transform of convolutions represents the product of Fourier transforms.
Somewhat surprisingly, for the study of convergence in $T_\infty$ an application 
of the $Q_h$-transform effectively replaces Fourier calculus. Indeed, the proof of 
Theorem 1.2 makes use of the Fourier analysis only in a minor way.

We will discuss the action of the $Q$-transform on single distributions
in Sections 2-5 and then turn to convolutions in Sections 6-7. The proof 
of Theorem 1.2 and Corollaries 1.3-1.4 are given in Section 8-9. 
The remaining Sections 10-14 illustrate these results for several
classes of probability distributions and specific examples.

\vskip2mm
{\sl Contents}

\vskip2mm
1. Introduction

2. Semigroup of shifted distributions. A basic identity

3. Moments of shifted distributions

4. Behaviour of eigenvalues near the critical zone

5. Maximum of shifted densities

6. Convolution of shifted densities

7. Local limit theorem for normalized shifted distributions

8. Proof of Theorem 1.2

9. Proof of Corollaries 1.3-1.4

10. Laplace transforms with separation property

11. Laplace transforms with periodic components

12. Characterization of periodicity in terms of Laplace transform

13. Periodic components via trigonometric series

14. Examples involving trigonometric polynomials

\vskip7mm
\section{{\bf Semigroup of Shifted Distributions. A Basic Identity}}
\setcounter{equation}{0}

\vskip2mm
\noindent
Let $X$ be a sub-Gaussian random vector in $\R^d$ with distribution $\mu$. 
Then, the Laplace transform, or the moment generating function
$$
(L \mu)(t) = L(t) = \E\,e^{\left<t,X\right>} = 
\int_{\R^d} e^{\left<t,x\right>}\,d\mu(x), \quad t \in \R^d,
$$
is finite and represents a $C^\infty$-smooth function on $\R^d$.
Correspondingly, the log-Laplace transform
$$
(K \mu)(t) = K(t) = \log L(t) = \log \E\,e^{\left<t,X\right>}
$$
is a convex, $C^\infty$-smooth function on $\R^d$.

The measures $\mu_h = Q_h \mu$ are defined in (1.4) by means 
of the Esscher transform. In the absolutely continuous case, Definition 1.5
is reduced to the following:

\vskip5mm
{\bf Definition 2.1.} If the random vector $X$ has density $p$, 
the measure $\mu_h$ has density which we denote similarly as
\be
Q_h p(x) = \frac{1}{L(h)}\,e^{\left<h,x\right>} p(x).
\en
In this case, let us also write $L\mu = Lp$ and $K\mu = Kp$.

\vskip5mm
We will call $\mu_h$ the shifted distribution of $X$ at step $h$
in order to emphasize the following 
important fact: For the standard normal density $\varphi(x)$, 
the shifted normal law has density $Q_h \varphi(x) = \varphi(x+h)$.

\vskip5mm
A remarkable property of the transform (1.4) is the semi-group property
$$
Q_{h_1}(Q_{h_2} \mu) = Q_{h_1 + h_2} \mu, \quad h_1,h_2 \in \R^d.
$$
Similarly, in the space of probability densities as in (2.1), we have
$$
Q_{h_1}(Q_{h_2} p) = Q_{h_1 + h_2} p.
$$

Let us also mention how this transform acts under rescaling, for simplicity
in the absolutely continuous case. Given $\lambda>0$,
the random vector $\lambda X$ has density 
$p_\lambda(x) = \lambda^{-d}\,p(x/\lambda)$ with Laplace transform
$(L p_\lambda)(t) = L(\lambda t)$. Hence
\be
Q_h p_\lambda(x) = \frac{1}{(L p_\lambda)(h)}\,e^{\left<h,x\right>} p_\lambda(x) = 
\frac{1}{\lambda^d}\, (Q_{\lambda h} p)\Big(\frac{x}{\lambda}\Big).
\en

The transform $Q_h$ is multiplicative with respect to convolutions.

\vskip5mm
{\bf Proposition 2.2.} {\sl If independent sub-Gaussian random
vectors in $\R^d$ have distributions $\mu_1,\dots,\mu_n$, then
for the convolution $\mu = \mu_1 * \dots * \mu_n$, we have
\be
Q_h \mu = Q_h \mu_1 * \dots * Q_h \mu_n.
\en
In particular, if $\mu_k$ have densities $p_k$, then
for the convolution $p = p_1 * \dots * p_n$, we have
$$
Q_h p = Q_h p_1 * \dots * Q_h p_n.
$$
}

{\bf Proof.} It is sufficient to compare the Laplace transforms  of both
sides in (2.3). The Laplace transform of $\mu$ is given by
$$
L \mu(t) = (L \mu_1)(t) \dots (L \mu_n)(t).
$$
Hence, the Laplace transform of $Q_h \mu$ is given by
\bee
(L Q_h \mu)(t) 
 & = & 
\int_{\R^d} e^{\left<t,x\right>}\,d Q_h \mu(x) \, = \,
\frac{1}{(L \mu)(t)} \int_{\R^d} e^{\left<t + h,x\right>} d\mu(x) \\
 & = &
\frac{(L \mu)(t+h)}{(L \mu)(t)} \, = \, 
\prod_{k=1}^n \frac{(L \mu_k)(t+h)}{(L \mu_k)(t)} \, = \, 
\prod_{k=1}^n (L Q_h \mu_k)(t).
\ene
\qed

\vskip2mm
If the random vector $X$ has density $p$, the formula (2.1) may be rewritten as
$$
p(x) = L(h) e^{-\left<x,h\right>}\, Q_h p(x) = 
e^{-\left<x,h\right> + K(h)}\, Q_h p(x),
$$
or
$$
\frac{p(x)}{\varphi(x)} = (2\pi)^{d/2}\,
e^{\frac{1}{2}\, |x-h|^2 - \frac{1}{2}\, |h|^2 + K(h)}\, Q_h p(x).
$$

Introduce the smooth function on $\R^d$
\be
(A p)(h) =  A(h) = \frac{1}{2}\, |h|^2 - K(h), 
\en
It allows one to reformulate the property $L(h) \leq e^{\frac{1}{2}\,|h|^2}$ 
as $A(h) \geq 0$ for all $h \in \R^d$. Note that this is equivalent to the strict 
sub-Gaussianity of $X$, when this random vector has mean zero and
identity covariance matrix. Thus, we have:

\vskip5mm
{\bf Proposition 2.3.} {\sl Given a sub-Gaussian random vector $X$ in $\R^d$ with
density $p$,
\be
\frac{p(x)}{\varphi(x)} = (2\pi)^{d/2}\,e^{\frac{1}{2}\, |x-h|^2 - A(h)}\, Q_h p(x),
\quad x,h \in \R^d,
\en
where $A$ is the associated function to $p$.
}

\vskip4mm
This is a basic identity which will be used with $h=x$ to bound from above 
the ratio on the left for densities of normalized sums of independent copies of $X$.

Since the function $A(h)$ plays an essential role in the representation (2.5), 
a number of its properties will be important in the sequel. Some of them may be explored
in the general situation where $X$ does not need to have a density.
We denote by $A'(h) = \nabla A(h)$ the gradient and by $A''(h)$ the Hessian of $A$ at the point $h$.

\vskip5mm
{\bf Proposition 2.4.} {\sl Let $X$ be a sub-Gaussian random vector in $\R^d$ 
with Laplace transform satisfying $A(h) \geq 0$ for all $h \in \R^d$. 
Then, for all $h \in \R^d$,
\be
|A'(h)|^2 \leq 2 A(h), 
\en
\be
A''(h) \leq I_d,
\en
\be
A(h) = 0 \, \Rightarrow \, A'(h) = 0 \ \, {\sl and} \ \, A''(h) \geq 0.
\en
}

{\bf Proof.} 
The inequalities in (2.7)-(2.8) are understood in matrix sense.
The implication in (2.8) is obvious. Also, since the function $K$ 
is convex, the assertion (2.7) follows from the definition (2.4). 

To prove (2.6), let us write
the Taylor integral formula for $A$ at the point $h$ up to the quadratic term
\be
A(h+x) = A(h) + \left<A'(h),x\right> + \int_0^1 (1-s) \left<A''(h+sx)x,x\right> ds.
\en
By (2.7), for all $x \in \R^d$,
$$
0 \leq A(h) + \left<A'(h),x\right> + \frac{1}{2}\,|x|^2.
$$
Minimizing the right-hand side over all $x$, we arrive at (2.6).
\qed

\vskip3mm
Note that an application of (2.6)-(2.7) in the Taylor formula (2.9) also implies that
\be
|A(h+x) - A(h)| \leq \sqrt{2A(h)}\,|x| + \frac{1}{2}\,|x|^2, \quad x,h \in \R^d.
\en

To explore some other properties of the function $A$, we need to look
at the moments of shifted distributions.

\vskip7mm
\section{{\bf Moments of Shifted Distributions}}
\setcounter{equation}{0}

\vskip2mm
\noindent
For a sub-Gaussian random vector $X$ in $\R^d$ with distribution $\mu$, denote 
by $X(h)$ a random vector with distribution $\mu_h = Q_h \mu$ ($h \in \R^d$). 
It is sub-Gaussian, and its Laplace and log-Laplace transforms are given by
\begin{eqnarray}
L_h(t) 
 & = & 
\E\,e^{\left<t,X(h)\right>} = \frac{L(t+h)}{L(h)}, \nonumber \\
K_h(t) 
 & = & 
\log L_h(t) = K(t+h) - K(h).
\end{eqnarray}

This random vector has mean
$$
m_h \, = \, \E X(h) \, = \, \frac{L'(h)}{L(h)} \, = \, K'(h),
$$
where we recall that $L'(h) = \nabla L(h)$ and $K'(h) = \nabla K(h)$ denote 
the gradients of $L$ and $K$ respectively. In addition, since for all 
$t,h \in \R^d$,
$$
\Var\big(\left<t,X(h)\right>\big) \, = \, 
\frac{d^2}{d\lambda^2}\,K_h(\lambda t)\Big|_{\lambda = 0} \, = \,
\left<K''(h) t,t\right>, 
$$
the matrix of second order partial derivatives (Hessian)
$$
R_h = K''(h), \quad h \in \R^d,
$$ 
represents the covariance matrix of $X(h)$. This shows that necessarily $K''(h)$ 
is positive semi-definite. Moreover, if $X$ has a density, it is strictly positive 
definite, since otherwise the random vector $X(h)$ may not have a density. 

The centered random vector
\be
\bar X(h) = X(h) - m_h
\en
has mean zero and covariance matrix $R_h$. Involving the function $A$ and 
assuming that it is non-negative, the log-Laplace transform is given by and satisfies
\bee
\bar K_h(t) \, = \,
\log \E\,e^{\left<t,\bar X(h)\right>} 
 & = &
K(t+h) - K(h) - \left<t,m_h\right> \\
 & = &
\Big(\frac{1}{2}\,|t+h|^2 - A(t+h)\Big) - \Big(\frac{1}{2}\,|h|^2 - A(h)\Big)  - 
\left<t,K'(h)\right> \\
 & \leq &
\frac{1}{2}\,|t+h|^2 - \Big(\frac{1}{2}\,|h|^2 - A(h)\Big)  - 
\left<t,K'(h)\right> \\
 & = &
\frac{1}{2}\,|t|^2 + A(h) + \left<t,A'(h)\right>.
\ene
Applying (2.6), we get
\bee
\bar K_h(t)
 & \leq &
\frac{1}{2}\,|t|^2 + A(h) + |t|\,|A'(h)| \\
 & \leq &
\frac{1}{2}\,|t|^2 + A(h) + |t| \sqrt{2 A(h)} \, = \,
\frac{1}{2}\,\Big(|t| + \sqrt{2 A(h)}\Big)^2.
\ene 
Thus, we obtain an upper bound for the Laplace transform.

\vskip5mm
{\bf Proposition 3.1.} {\sl Let $X$ be a sub-Gaussian random vector in $\R^d$ 
with Laplace transform such that $A(h) \geq 0$ for all $h \in \R^d$. 
Then
\be
\E\,e^{\left<t,\bar X(h)\right>} \leq
\exp\Big\{\frac{1}{2}\Big(|t| + \sqrt{2 A(h)}\Big)^2\Big\}, \quad t,h \in \R^d.
\en
}

This relation shows in a quantitative form the sub-Gaussian character of $\bar X$.
Indeed, to simplify, one may use $(a+b)^2 \leq 2 a^2 + 2b^2$ ($a,b \in \R$), 
so that (3.3) yields
$$
\E\,e^{\left<t,\bar X(h)\right>} \leq \exp\big\{|t|^2 + 2A(h)\big\}.
$$
Let us apply this with $t = \frac{1}{2}\,\xi \theta$, $|\theta| = 1$, where $\xi$ is 
a standard normal random variable independent of $X$. Taking the expectation 
with respect to $\xi$ and using $\E\,e^{\frac{1}{4}\,\xi^2} = \sqrt{2}$, we then get
\be
\E\,e^{\frac{1}{4} \left<\theta,\bar X(h)\right>^2} \leq
\sqrt{2}\,e^{2 A(h)}.
\en

By Jensen's inequality, the left expectation is greater than
$\exp\big\{\frac{1}{4}\, \E \left<\theta,\bar X(h)\right>^2\big\}$ which leads to
$$
\frac{1}{4}\left<R_h \theta,\theta\right> \leq 2 A(h) + \log\sqrt{2}.
$$
In particular, one may apply this inequality to orthonormal eigenvectors 
$\theta$ of $R_h$.

\vskip5mm
{\bf Corollary 3.2.} {\sl Under the conditions of Proposition $3.1$, 
the eigenvalues $\lambda_j(h)$ of $R_h$ satisfy, for any $h \in \R^d$, 
\be
\max_{1 \leq j \leq d} \lambda_j(h) \, \leq \, 8A(h) + 2\log 2.
\en
}

\vskip2mm
Starting from (3.4), one may similarly estimate higher order moments
of linear functionals of $\bar X(h)$. The following bound will be needed
with $q=3$.

\vskip5mm
{\bf Corollary 3.3.} {\sl Under the conditions of Proposition $3.1$, 
up to some absolute constant $C$,
\be
\Big(\E\,\big|\left<\theta,\bar X(h)\right>\big|^q\Big)^{2/q} \leq 
C\big(q + A(h)\big), \quad q \geq 1,
\en
}

{\bf Proof.} One may assume that $q \geq 2$. By (3.4), for the random variable 
$$
\eta = \frac{1}{4} \left<\theta,\bar X(h)\right>^2 - 2 A(h) - \log \sqrt{2}
$$
we have $\E\,e^\eta \leq 1$, implying $\P\{\eta \geq x\} \leq e^{-x}$
for all $x \geq 0$. Hence, for any $r \geq 1$,
$$
\E\,(\eta^+)^r = r \int_0^\infty x^{r-1}\,\P\{\eta \geq x\}\,dx \leq
\Gamma(r+1),
$$
where $\eta^+ = \max(\eta,0)$. Using 
$\left<\theta,\bar X(h)\right>^2 \leq 4 \eta^+ + 8 A(h) + 2\log 2$,
it follows that
$$
\Big(\E\,\big|\left<\theta,\bar X(h)\right>\big|^{2r}\Big)^{1/r} \leq 
4\,\Gamma(r+1)^{1/r}+ 8 A(h) + 2\log 2.
$$
It remains to apply this inequality with $r = q/2$.
\qed

\vskip7mm
\section{{\bf Behaviour of Eigenvalues near the Critical Zone}}
\setcounter{equation}{0}

\vskip2mm
\noindent
Another application of the sub-Gaussian bound (3.4), more precisely -- of
the inequality (3.6) with $q=3$, concerns the behaviour
of the Hessian $A''(h)$ when $A(h)$ is small (which we shall call the
critical zone in view of the inequality (6.3) from Proposition 6.3 below). 
The following statement complements the property (2.8) from 
Proposition 2.4 which asserts that $A''(h) \geq 0$ as long as $A(h) = 0$.

\vskip5mm
{\bf Proposition 4.1.} {\sl Let $X$ be a sub-Gaussian random vector in $\R^d$ 
with Laplace transform such that $A(h) \geq 0$ for all $h \in \R^d$. 
If $0 \leq A(h) \leq 1$, then
\be
\inf_{|\theta| = 1} \left<A''(h) \theta,\theta\right> \geq -C_d A(h)^{1/4}
\en
with some constant $C_d>0$ depending on $d$ only.
}

\vskip4mm
One consequence of (4.1), which will be needed for the characterization 
of the CLT with respect to the R\'enyi divergence of infinite order is that
$$
\liminf_{A(h) \rightarrow 0} \ \inf_{|\theta| = 1} \left<A''(h) \theta,\theta\right> 
 \geq 0,
$$
or equivalently, since $A'' = I_d - K''$,
$$
\limsup_{A(h) \rightarrow 0} \ \sup_{|\theta| = 1} \left<K''(h) \theta,\theta\right> 
 \leq 1.
$$
In terms of the eigenvalues $\lambda_j(h)$ of the covariance matrix 
$R_h$ of the random vector $X(h)$ this may also be restated as:

\vskip4mm
{\bf Corollary 4.2.} {\sl Let $X$ be a sub-Gaussian random vector in $\R^d$ 
with Laplace transform such that $A(h) \geq 0$ for all $h \in \R^d$. Then
\be
\limsup_{A(h) \rightarrow 0} \ \max_{1 \leq j \leq d} \lambda_j(h)
 \leq 1.
\en
In particular,
$$
\limsup_{A(h) \rightarrow 0} \ {\rm det} K''(h) \leq 1.
$$
}

{\bf Proof.} We keep notations from the previous sections.
Recall that
$$
\bar X(h) = X(h) - \E X(h) = (\xi_1,\dots,\xi_d), \quad
h = (h_1, \dots,h_d) \in \R^d.
$$

The following elementary identity is needed for the 3rd order partial
derivatives of the log-Laplace transform with respect to the variables
$h_i$, $h_j$, $h_k$:
$$
\partial_{ijk} K(h) = \E\, \xi_i \xi_j \xi_k = -\partial_{ijk} A(h), 
\quad 1 \leq i,j,k \leq d.
$$
By (3.6) with $q=3$, 
$$
\big(\E\,|\xi_i|^3\big)^{1/3} \leq C\,(1+A(h))^{1/2},
$$
where $C$ is an absolute constant. Hence, by H\"older's inequality, 
for any $h \in \R^d$,
\be
|\partial_{ijk} A(h)| \leq C^3\,(1+A(h))^{3/2}.
\en

Let $0 < A(h) \leq 1$. Applying the inequality (2.10),
we get that, whenever $|x| \leq 1$,
$$
|A(h + x) - A(h)| \leq 
\sqrt{2 A(h)}\,|x| + \frac{1}{2}\,|x|^2 \, < \, 2.
$$
Hence $A(h + x) < 3$, so that by (4.3), 
\be
|\partial_{ijk} A(h + x)| \leq C
\en
with some absolute constant $C>0$.

We now apply the multidimensional integral Taylor formula up
to cubic terms which indicates that
\begin{eqnarray}
A(h+x) 
 & = &
A(h) + \left<A'(h),x\right> + \frac{1}{2}\left<A''(h)x,x\right> \nonumber \\
 & & + \
\sum_{|\beta| = 3} \frac{3}{\beta!}\,x^\beta \int_0^1 
(1-s)^2\, D^\beta A(h + sx)\,ds.
\end{eqnarray}
Here we use the standard notation for the partial derivative 
$$
D^\beta A = 
\frac{\partial^{|\beta|} A}{\partial h_1^{\beta_1} \dots \partial h_d^{\beta_d}},
\quad \beta = (\beta_1,\dots,\beta_d),
$$ 
where $\beta$ is a multi-index of length 
$|\beta| = \beta_1 + \dots + \beta_d$ (with integers $\beta_i \geq 0$), 
and where $x = x_1^{\beta_1} \dots x_d^{\beta_d}$ for 
$x = (x_1,\dots,x_d) \in \R^d$. 

If $|x| \leq 1$, then, by (4.4), all these partial derivatives for $|\beta| = 3$ 
do not exceed $C$ in absolute value. Since also $|x^\beta| \leq |x|^{|\beta|}$, 
we obtain from (4.5) that
\be
A(h+x) \, \leq \,
A(h) + \sqrt{2A(h)}\,|x| + \frac{1}{2}\left<A''(h)x,x\right> + C_d\,|x|^3
\en
with
$$
C_d = C \sum_{|\beta| = 3} \frac{1}{\beta!}.
$$
Here we also used (2.6) in order to bound the linear term in (4.5).
Let us choose $x = r\theta$ with $0 < r \leq 1$, $|\theta|=1$, 
and use the assumption $A(h+x) \geq 0$. Then (4.6) yields
$$
\left<A''(h)\theta,\theta\right> \geq 
-\frac{2}{r^2}\,\Big[A(h) + \sqrt{2 A(h)}\,r + C_d r^3\Big].
$$
In particular, the choice $r  = A(h)^{1/4}$ leads to (4.1).
\qed

\vskip7mm
\section{{\bf Maximum of Shifted Densities}}
\setcounter{equation}{0}

\vskip2mm
\noindent
In order to bound the last term $Q_h p(x)$ in the basic identity (2.5), 
suppose that the distribution of $X$ has a finite R\'enyi distance of infinite order 
to the standard normal law. This means that $X$ has a density $p$ which admits 
a pointwise upper bound
\be
p(x) \leq c \varphi(x), \quad x \in \R^d \ ({\rm a.e.})
\en
with optimal value $c = 1 + T_\infty(p||\varphi)$. 
In that case, one may control the maximum 
$$
M(Q_h p) = {\rm ess\,sup}_x\, Q_h p(x)
$$
as follows. By (5.1), for almost all $x \in \R^d$,
\bee
Q_h p(x)
 & = &
\frac{1}{L(h)}\,e^{\left<x,h\right>} p(x) \\
 & \leq &
\frac{c\, e^{\left<x,h\right> - \frac{1}{2}\, |x|^2}}{L(h)\,(2\pi)^{d/2}}\,  \leq \,
\frac{c\, e^{\frac{1}{2}\,|h|^2}}{L(h)\, (2\pi)^{d/2}} \, = \,
\frac{c}{(2\pi)^{d/2}}\, e^{A(h)},
\ene
where $L$ is the Laplace transform of the distribution of $X$ and 
$A(h) = \frac{1}{2}\,|h|^2 - K(h)$, $K(h) = \log L(h)$. Thus, we arrive at the
following elementary relation.

\vskip5mm
{\bf Proposition 5.1.} {\sl Let $X$ be a sub-Gaussian random vector in $\R^d$ 
with density $p$ such that $c = 1 + T_\infty(p||\varphi)$ is finite. Then, 
for all $h \in \R^d$,
\be
M(Q_h p) \leq \frac{c}{(2\pi)^{d/2}}\, e^{A(h)}.
\en
}

This inequality may be used to bound the eigenvalues of the covariance 
matrix $R_h$ in terms of $A(h)$ from below. 
This complements the upper bound of Corollary 3.2. 

In the sequel, put
$$
\sigma_h = ({\rm det}(R_h))^{\frac{1}{2d}} = 
({\rm det}(K''(h))^{\frac{1}{2d}}, \quad h \in \R^d,
$$
which is everywhere positive. In dimension $d=1$, this quantity represents
the standard deviation of the random variable $X(h)$.

\vskip4mm
{\bf Proposition 5.2.} {\sl Under the condition $(5.1)$, for all $h \in \R^d$,
\be
\sigma_h^d \geq \frac{1}{c}\, e^{-A(h) - d/2}.
\en
}

{\bf Proof.} The argument employs the following known relation between the
covariance matrix and maximum of density (whose proof we include for
completeness at the end of this section): Given a random vector $\xi$ in $\R^d$ 
with finite second moment and finite $M = M(\xi)$, we have
\be
\big(M^2\,{\rm det}(R)\big)^{\frac{1}{d}} \geq \frac{1}{2\pi e},
\en
where $R$ is the covariance matrix of $\xi$. Applying this inequality to 
the random vector $\xi = X(h)$ with its covariance matrix $R = R_h$ and 
using (5.2), we get
$$
\frac{1}{(2\pi e)^{d/2}} \leq M(X(h)) \,\sigma_h^d  \leq 
\frac{c\,\sigma_h^d}{(2\pi)^{d/2}}\, e^{A(h)},
$$
from which (5.3) follows immediately.
\qed

\vskip4mm
{\bf Corollary 5.3.} {\sl Let $X$ be a sub-Gaussian random vector in $\R^d$ 
with Laplace transform such that $A(h) \geq 0$ for all $h \in \R^d$ and with
finite $c = 1 + T_\infty(p||\varphi)$. Then,
the eigenvalues $\lambda_j(h)$ of $R_h$ satisfy, for any $h \in \R^d$, 
\begin{eqnarray}
\frac{1}{c^2}\, \frac{e^{-2A(h) - d}\,}{(8A(h) + 2)^{d-1}}
 & \leq &
\min_{1 \leq j \leq d} \lambda_j(h) \nonumber \\
 & \leq & 
\max_{1 \leq j \leq d} \lambda_j(h) \, \leq \, 8A(h) + 2.
\end{eqnarray}
}

\vskip2mm
{\bf Proof.} Put $\lambda_j = \lambda_j(h)$ and $\alpha = A(h)$. 
The upper bound in (5.5) is provided in Lemma 3.2. On the other hand, by (5.3),
\be
\sigma_h^{2d} = {\rm det}(R_h) = \lambda_1 \dots \lambda_d \geq  
\frac{1}{c^2}\, e^{-2\alpha - d}.
\en
Therefore, by the upper bound in (5.5),
$$
\lambda_1 \dots \lambda_d \leq \min_j \lambda_j\, (\max_j \lambda_j)^{d-1}
 \leq \min_j \lambda_j\,(8\alpha + 2)^{d-1},
$$
and (5.6) yields the lower bound in (5.5).
\qed

\vskip4mm
{\bf Proof of (5.4).} We follow a simple information-theoretic approach
proposed in \cite{B-M}. Introduce the entropy functional
\be
h(p) = -\int_{\R^d} p(x) \log p(x)\,dx.
\en
It is well defined for absolutely continuous distributions with finite second moment 
and is maximized for the normal distribution when the covariance matrix is fixed.
Indeed, without loss of generality, let $\xi$ have mean zero. If 
$\zeta$ has a normal density $q$ on $\R^d$ with mean zero and the same 
covariance matrix $R$, then 
$$
\E \left<R^{-1} \xi,\xi\right> = \E \left<R^{-1} \zeta,\zeta\right> = d,
$$
implying
$$
h(q) - h(p) = \int_{\R^d} p(x) \log \frac{p(x)}{q(x)}\,dx.
$$
Here the right-hand side defines the relative entropy $D(p||q)$ of $p$ 
with respect to $q$ (the Kullback-Leibler distance), which is non-negative, 
by Jensen's inequality.

Now, from one hand,
$$
h(q) = \frac{d}{2}\,\log(C\sigma^2), \quad C = 2\pi e, \ \
\sigma = ({\rm det}(R))^{\frac{1}{2d}}.
$$
On the other hand, $h(p) \geq -\log M$, which follows from (5.7)
using $\log p(x) \leq \log M$ (a.e.) Hence
$$
D(p||q) \leq \frac{d}{2}\,\log(C\sigma^2) + \log M,
$$
that is,
$$
M^2\,{\rm det}(R) \geq \frac{1}{(2\pi e)^d}\,
e^{2 D(q||p)}.
$$
This is a sharpened form of (5.4).
\qed

\vskip7mm
\section{{\bf Convolutions of Shifted Distributions}}
\setcounter{equation}{0}

\vskip2mm
\noindent
We are now prepared to apply these results to the normalized sums
$$
Z_n = \frac{X_1 + \dots + X_n}{\sqrt{n}}
$$
of independent copies of a sub-Gaussian random vector $X$ in $\R^d$ 
with density $p$. Note that in terms of $L = L p$, $K = K p$ and $A = A p$, 
for the density $p_n$ of $Z_n$ we have
\bee
(L p_n)(t)
 & = &
L(t/\sqrt{n})^n = e^{n K(t/\sqrt{n})}, \\
(K p_n)(t) 
 & = &
n K(t/\sqrt{n}), \\
(A p_n)(h_n) 
 & = & 
\frac{1}{2}\, |h_n|^2 - (Kp_n)(h_n) \, = \, \frac{n}{2}\, |h|^2 - n K(h) \, = \, n A(h),
\ene
where $h_n = h\sqrt{n}$. Hence, the basic identity (2.5) yields a similar formula.

\vskip5mm
{\bf Proposition 6.1.} {\sl Putting $x_n = x\sqrt{n}$, $h_n = h\sqrt{n}$
$(x,h \in \R^d)$, we have
\be
\frac{p_n(x\sqrt{n})}{\varphi(x\sqrt{n})} = 
(2\pi)^{d/2}\,e^{\frac{n}{2}\, |x - h|^2 - n A(h)}\, Q_{h_n} p_n(x_n).
\en
}

This equality becomes useful, if we are able to bound the factor $Q_{h_n} p_n(x_n)$ 
uniformly over all $x$ for a fixed value of $h$ as stated in the following:

\vskip5mm
{\bf Corollary 6.2.} {\sl For all $x,h \in \R^d$,
\be
\frac{p_n(x\sqrt{n})}{\varphi(x\sqrt{n})} \leq 
(2\pi)^{d/2}\,e^{\frac{n}{2}\, |x - h|^2 - n A(h)}\, M(Q_{h\sqrt{n}}\, p_n).
\en
}

\vskip2mm
A general upper bound on the $M$-functional was given in Proposition 5.1. 
However, it is useless to apply this bound directly in (6.2) to the densities $p_n$, 
since then the right-hand side of (5.2) will contain the parameter 
$c_n = 1 + T_\infty(p_n||\varphi)$. Instead, we use a semi-additive property 
of the maximum-of-density functional, which indicates that
$$
M(X_1 + \dots + X_n)^{-\frac{2}{d}} \geq 
\frac{1}{e} \sum_{k=1}^n M(X_k)^{-\frac{2}{d}}
$$
for all independent random vectors $X_k$ in $\R^d$ having bounded 
densities, cf. \cite{B-C}. If all $X_k$ 
are identically distributed and have a density $p$, this relation yields
$$
M(p^{*n}) \leq \Big(\frac{e}{n}\Big)^{d/2} M(p)
$$
for the convolution $n$-th power of $p$. Applying the multiplicativity
property of the Esscher transform (Proposition 2.2) together with (5.2), 
we then have
$$
M(Q_h p^{*n}) \leq \Big(\frac{e}{n}\Big)^{d/2} M(Q_h p) \leq 
c\,\Big(\frac{e}{2 \pi n}\Big)^{d/2}\, e^{A(h)},
$$
where we recall that $c = 1 + T_\infty(p||\varphi)$.
Now, since $p^{*n}(x) = \frac{1}{\lambda}\,p_n(\frac{x}{\lambda})$ with
$\lambda = \sqrt{n}$, one may apply the scaling identity (2.2):
$$
M(Q_h p^{*n}) = \frac{1}{n^{d/2}}\, M(Q_{h\sqrt{n}}\, p_n).
$$
Hence
$$
M(Q_{h\sqrt{n}}\, p_n) \leq c\, \Big(\frac{e}{2\pi}\Big)^{d/2}\, e^{A(h)}.
$$
One can now return to Corollary 6.2 and apply the above bound to get that
$$
\frac{p_n(x\sqrt{n})}{\varphi(x\sqrt{n})} \leq 
c\,e^{d/2}\,e^{\frac{n}{2}\, |x - h|^2 - (n-1) A(h)}.
$$
In particular, with $h=x$ this yields:

\vskip5mm
{\bf Proposition 6.3.} {\sl If the density $p$ has finite R\'enyi distance 
of infinite order to the standard normal law, then, for almost all $x \in \R^d$,
\be
\frac{p_n(x\sqrt{n})}{\varphi(x\sqrt{n})} \leq 
c\,e^{d/2}\,e^{- (n-1) A(x)},
\en
where $c = 1 + T_\infty(p||\varphi)$.
}

\vskip5mm
{\bf Corollary 6.4.} {\sl If additionally $X$ has identity covariance
matrix and is strictly sub-Gaussian, then
$$
T_\infty(p_n||\varphi) \leq e^{d/2}\,(1 + T_\infty(p||\varphi)) - 1.
$$
}

\vskip2mm
Thus, the finiteness of the Tsallis distance $T_\infty(p||\varphi)$ for a strictly 
sub-Gaussian random vector $X$ with density $p$ ensures the boundedness
of $T_\infty(p_n||\varphi)$ for all normalized sums $Z_n$.

If $A(x)$ is bounded away from zero, the inequality (6.3) shows that
$p_n(x\sqrt{n})/\varphi(x\sqrt{n})$ is exponentially small for growing $n$.
In particular, this holds for any non-normal random vector $X$ satisfying 
the separation property (1.3). Then we immediately obtain:

\vskip5mm
{\bf Corollary 6.5.} {\sl Suppose that $X$ has a density $p$ with finite 
$T_\infty(p||\varphi)$. Under the condition $(1.3)$, for any $\tau_0 > 0$, 
there exist $A > 0$ and $\delta \in (0,1)$ such that the densities $p_n$ 
of $Z_n$ satisfy
\be
p_n(x) \leq A \delta^n \varphi(x), \quad |x| \geq \tau_0\sqrt{n}.
\en
}

In particular,
\be
\liminf_{n \rightarrow \infty} \, \sup_{x \in \R^d} \,
\frac{|p_n(x) - \varphi(x)|}{\varphi(x)} \geq 1.
\en
Therefore, one can not hope to strengthen the Tsallis distance by
introducing a modulus sign in the definition of the distance.

Although (1.3) does not hold in general for strictly sub-Gaussian distributions, 
Proposition 6.3 may be applied outside the set of points where $A(x)$ is 
bounded away from zero. If $A(x)$ is close to zero, we say that the point 
$x$ belongs to the critical zone (a precise definition will be given in Section 8).
In this case, we need to return to the basic representation of Proposition 6.1
and  study the last term $Q_{h_n} p_n(x_n)$. This requires to apply a variant
of the local limit theorem, using the property that the density $Q_{h_n} p_n$ has
a convolution structure.

\vskip7mm
\section{{\bf Local Limit Theorem for Normalized Shifted Distributions}}
\setcounter{equation}{0}

\vskip2mm
\noindent
Keeping the notations and assumptions from the previous section,
recall that $R_h = K''(h)$ represents the covariance matrix of $X(h)$.
Since it is symmetric and strictly positive definite, one may consider 
the normalized random vectors
\bee
\widehat X(h) 
 & = &
R_h^{-1/2} (X(h) - \E X(x)) \nonumber \\
 & = & 
R_h^{-1/2} (X(h) - m_h), \qquad h \in \R^d,
\ene
which have mean zero and identity covariance matrix.
We will have to consider convolution powers of distributions of $\widehat X(h)$
by means of a multidimensional local limit theorem for the points $h$ where 
the value 
$$
A(h) = \frac{1}{2}\,|h|^2 - K(h) = \frac{1}{2}\,|h|^2 - \log \E\,e^{\left<X,h\right>}
$$  
is small. More precisely, here we prove the following 
refinement of the representation (6.1).
Consider the vector-function
\begin{eqnarray}
v_x 
 & = &
R_x^{-1/2}(x - m_x) \nonumber \\
 & = &
R_x^{-1/2}(x - K'(x)) \, = \, R_x^{-1/2} A'(x), \quad x \in \R^d,
\end{eqnarray}
and recall that
$$
\sigma_x = ({\rm det} R_x)^{\frac{1}{2d}} = 
({\rm det} K''(x))^{\frac{1}{2d}}.
$$

\vskip5mm
{\bf Proposition 7.1.} {\sl If the Laplace transform of a sub-Gaussian random 
vector $X$ in $\R^d$ with finite $c = 1 + T_\infty(p||\varphi)$
is such that $A(h) \geq 0$ for all $h \in \R^d$, then for all $x \in \R^d$, $n \geq 6$,
we have
\be
\frac{p_n(x\sqrt{n})}{\varphi(x\sqrt{n})} = 
\frac{1}{\sigma_x^d}\,\exp\big\{-nA(x)-n\,|v_x|^2/2\big\} + \frac{B^d c^5}{\sqrt{n}},
\en
where $B = B_n(x)$ is bounded by an absolute constant.
}

\vskip5mm
One should note that the term $|v_x|^2$ appearing on the right-hand side 
of (7.2) must be small for small values of $A(x)$.
Indeed, according to Corollary 5.3, the minimal eigenvalue 
$\lambda = \min_{1 \leq j \leq d} \lambda_j(x)$ of the matrix $R_x$ admits
a lower bound
\be
\lambda \geq
\frac{1}{c^2\,(8A(x) + 2)^{d-1}} \, e^{-2A(x) - d}.
\en
As a consequence, applying Proposition 2.4 in (7.1), we have
\begin{eqnarray}
|v_x|^2 
 & \leq &
\frac{1}{\lambda}\,|A'(x)|^2 \nonumber \\
 & \leq &
2c^2\,(8A(x) + 2)^{d-1} \, A(x)\, e^{2A(x) + d} \, \leq \,
C^d c^2 A(x)
\end{eqnarray}
for some absolute constant $C>0$, where we assumed that $A(x) \leq 1$
in the last step.

For the derivation of (7.2), we employ a general local limit theorem 
for densities on $\R^d$ with a quantitative error term, which was 
recently derived in \cite{B-G2}.

\vskip5mm
{\bf Lemma 7.2.} {\sl Let $(\xi_k)_{k \geq 1}$ be independent copies of 
a random vector $\xi$ in $\R^d$ with mean zero, identity covariance matrix and 
finite third absolute moment. Assuming that $\xi$ has a bounded density,
the densities $q_n$ of the normalized sums
$Z_n = (\xi_1 + \dots + \xi_n)/\sqrt{n}$ satisfy 
$$
\sup_x |q_n(x) - \varphi(x)| \leq C^d\,\frac{1}{\sqrt{n}}\,M(\xi)^2\,\E\,|\xi|^3, \quad
x \in \R^d,
$$
with some absolute constant $C>0$.
}

\vskip5mm
{\bf Proof of Proposition 7.1.}
Consider the term $Q_{h_n} p_n$ in (6.1) with $h_n = h\sqrt{n}$. 
By Proposition 3.2, this density has a convolution structure.
It was also emphasized in (2.2) that, for any random vector $X$ with density $p = p_X$,
$$
Q_h p_{\lambda X}(x) = 
\frac{1}{\lambda^d}\,(Q_{\lambda h} p)\Big(\frac{x}{\lambda}\Big).
$$
Using this notation, we have $p_n = p_{S_n/\sqrt{n}}$ for the sum
$S_n = X_1 + \dots + X_n$. Hence with $\lambda = 1/\sqrt{n}$, 
$$
Q_{h_n} p_n(x) = n^{d/2}\, (Q_h p_{S_n})(x\sqrt{n}) = 
n^{d/2}\, (Q_h p) * \dots * (Q_h p) (x\sqrt{n}),
$$
where we applied Proposition 2.2 in the last step. Since $Q_h p$ serves
as a density of the random vector $X(h)$, $Q_{h_n} p_n(x)$ represents
the density for the normalized sum
$$
Z_{n,h} = \frac{X_1(h) + \dots + X_n(h)}{\sqrt{n}},
$$
assuming that $X_k(h)$ are independent. Introduce the normalized sums
$$
\widehat Z_{n,h} = \frac{\widehat X_1(h) + \dots + \widehat X_n(h)}{\sqrt{n}}
$$
for the shifted distributions, i.e. with
$X_k(h) = m_h + R_h^{1/2} \widehat X_k(h)$. Thus,
$$
Z_{n,h} = m_h \sqrt{n} + R_h^{1/2} \widehat Z_{n,h}.
$$
Denote by $\widehat p_{n,h}$ the density of $\widehat Z_{n,h}$. Then the density of
$Z_{n,h}$ is given by
$$
p_{n,h}(x) \, = \,
\frac{1}{\sigma_h^d}\,\widehat p_{n,h}\big(R_h^{-1/2} (x - m_h \sqrt{n})\big).
$$
At the points $x_n = x\sqrt{n}$ as in (6.1), we therefore obtain that
$$
Q_{h_n} p_n(x_n) = p_{n,h}(x_n) = 
\frac{1}{\sigma_h^d}\,\widehat p_{n,h}\big(\sqrt{n}\,R_h^{-1/2}(x - m_h)\big).
$$
Consequently, the equality (6.1) may be equivalently stated as
$$
\frac{p_n(x\sqrt{n})}{\varphi(x\sqrt{n})} = 
(2\pi)^{d/2}\,e^{\frac{n}{2}\, (x - h)^2 - n A(h)}\,
\frac{1}{\sigma_h^d}\,\widehat p_{n,h}\big(\sqrt{n}\,R_h^{-1/2}(x - m_h)\big).
$$
In particular, for $h = x$, we get
\be
\frac{p_n(x\sqrt{n})}{\varphi(x\sqrt{n})} = (2\pi)^{d/2}\,e^{- n A(x)}\,
\frac{1}{\sigma_x^d}\,\widehat p_{n,x}\big(v_x\sqrt{n}\big).
\en

We are in a position to apply Lemma 7.2 to the sequence 
$\xi_k = \widehat X_k(x)$ and write
\be
\widehat p_{n,x}(z) = \varphi(z) + B^d\,\frac{\beta_3(x) M_x^2}{\sqrt{n}},
\quad z \in \R^d,
\en
where
$$
\beta_3(x) = \E\,\big|\widehat X(x)\big|^3, \quad
M_x = M\big(\widehat X(x)\big),
$$
and where the quantity $B = B_n(z)$ is bounded by an absolute constant.
The latter maximum $M_x$ can be bounded by virtue of the upper bound (5.2):
$$
M(\widehat X(x)) = \sigma_x^d\, M(X(x)) = \sigma_x^d\, M(Q_x p)
\leq \frac{c\sigma_x^d}{(2\pi)^{d/2}}\, e^{A(x)}.
$$
In this case, (7.6) may be simplified with a new quantity $B = B_n(z)$ to
$$
\widehat p_{n,x}(z) = \varphi(z) + 
B^d c^2\,\frac{\beta_3(x) \sigma_x^{2d}}{\sqrt{n}}\, e^{2A(x)}.
$$
Inserting this in (7.5) with $z = v_x \sqrt{n}$, we arrive at
\be
\frac{p_n(x\sqrt{n})}{\varphi(x\sqrt{n})} = \frac{1}{\sigma_x^d}\,
e^{-nA(x) -n |v_x|^2/2} + 
B^d c^2\,\frac{\beta_3(x) \sigma_x^d}{\sqrt{n}}\, e^{-(n-2)A(x)},
\en
where $B = B_n(x)$ is bounded in absolute value by an absolute constant.

In order to estimate $\beta_3(x)$, first let us recall the bound
of Corollary 3.3 which implies
\be
\E\,|X(x) - m_x|^3 \leq Cd^{3/2}\,(1+A(x))^{3/2}
\en
with some absolute constant $C$. Also, by (7.3),
for any $w \in \R^d$,
$$
\big|R_x^{-1/2} w\big| \leq
c\,(8A(x) + 2)^{\frac{d-1}{2}} \, e^{A(x) + d/2}\,|w|.
$$
Applying this with $w = X(x) -m_x$ together with (7.8), we get
\bee
\beta_3(x) 
 & \leq &
c^3\,(8A(x) + 2)^{\frac{3(d-1)}{2}} \, e^{3A(x) + 3d/2}\cdot
Cd^{3/2}\,(1+A(x))^{3/2} \\
 & \leq &
C_1^d c^3\,e^{4A(x)}
\ene
for some absolute constant $C_1>0$. It remains to insert this bound
in (7.7) leading to
$$
\frac{p_n(x\sqrt{n})}{\varphi(x\sqrt{n})} = \frac{1}{\sigma_x^d}\,
e^{-nA(x) -n |v_x|^2/2} + 
B^d c^5\,\frac{\sigma_x^d}{\sqrt{n}}\, e^{-(n-6)A(x)}.
$$
\qed

\vskip7mm
\section{{\bf Proof of Theorems 1.2}}
\setcounter{equation}{0}

\vskip2mm
\noindent
Recall that the assumptions 1)-2) stated before Theorem 1.2 are necessary
for the convergence $T_\infty(p_n||\varphi) \rightarrow 0$ as 
$n \rightarrow \infty$. For simplicity, we assume that $n_0 = 1$, that is, 
$X$ is a strictly subgaussian random vector with mean zero, identity
covariance matrix, and finite constant $c = 1 + T_\infty(p||\varphi)$. 
In particular, the function
$$
A(x) = \frac{1}{2}\, |x|^2 - K(x)
$$
is non-negative on the whole space $\R^d$. 

According to Corollary 4.2, the function 
$\sigma_x = ({\rm det}\, K''(x))^{\frac{1}{2d}}$ satisfies
\be
\limsup_{A(x) \rightarrow 0}\, \sigma_x^{2d} \leq 1.
\en
First we show that the convergence 
$T_\infty(p_n||\varphi) \rightarrow 0$ is equivalent to 
\be
\lim_{A(x) \rightarrow 0}  \sigma_x^{2d} = 
\lim_{A(x) \rightarrow 0} {\rm det}\, K''(x) = 1,
\en
which is a compact version of the conditions $a')-b')$ mentioned after Theorem 1.2. 

Introduce the critical zones 
$$
A_n(a) = \Big\{x \in \R^d: A(x) \leq \frac{a}{n-1}\Big\}, \quad a>0.
$$

\vskip2mm
{\bf Sufficiency part.} Putting
$a = \log(1/\ep)$, $\ep \in (0,1)$, the upper bound (6.3) yields
$$
\sup_{x \notin A_n(a)}\, \frac{p_n(x\sqrt{n})}{\varphi(x\sqrt{n})} \leq ce^{d/2}\,\ep.
$$
As for the critical zone, the equality (7.2) is applicable for $n \geq 6$ and implies
$$
\sup_{x \in A_n(a)}\, \frac{p_n(x\sqrt{n})}{\varphi(x\sqrt{n})} \leq 
\sup_{x \in A_n(a)}\, \frac{1}{\sigma_x^d} + O\Big(\frac{1}{\sqrt{n}}\Big).
$$
Both estimates can be combined to give
$$
1 + T_\infty(p_n||\varphi) \leq \sup_{x \in A_n(a)}\, \frac{1}{\sigma_x^d} + 
ce^{d/2}\,\ep + O\Big(\frac{1}{\sqrt{n}}\Big).
$$
Thus, a sufficient condition for the convergence $T_\infty(p_n||\varphi) \rightarrow 0$ 
as $n \rightarrow \infty$ is that, for any $\ep \in (0,1)$,
$$
\limsup_{n \rightarrow \infty}\, \sup_{x \in A_n(\log(1/\ep))}\, 
\sigma_x^{-d} \leq 1.
$$
Since $A(x) = O(1/n)$ on every set $A_n(a)$, the above may be written as
the condition
\be
\liminf_{A(x) \rightarrow 0} \, \sigma_x^{2d} \geq 1,
\en
which is equivalent to (8.2) in view of (8.1).

\vskip2mm
{\bf Necessity part.} 
To see that the condition (8.2) is also necessary for the convergence in $T_\infty$,
let us return to the representation (7.2). Assuming that 
$T_\infty(p_n||\varphi) \rightarrow 0$, it implies that, for any $a>0$,
\be
\limsup_{n \rightarrow \infty}\, \sup_{x \in A_n(a)}\, \frac{1}{\sigma_x^d}
\exp\Big\{- n \Big(A(x) + \frac{1}{2}\,|v_x|^2\Big)\Big\} \leq 1.
\en
Recall that $|v_x|^2 \, \leq \, C^d c^2 A(x)$ whenever $n \geq a+1$,
as explained in (7.4). Since  $nA(x) \leq 2a$ on the set $A_n(a)$, 
it follows that 
$$
A(x) + \frac{1}{2}\,|v_x|^2 \leq \beta A(x) \leq \frac{2\beta}{n}\,a,
$$
where $\beta$ depends on $c$ and $d$ only. Thus, (8.4) implies that
$$
\limsup_{n \rightarrow \infty}\, \sup_{x \in A_n(a)}\, \frac{1}{\sigma_x^{2d}}
\leq e^{4\beta a}, \quad 0 < a \leq 1.
$$
Therefore, for all $n \geq n(a)$,
$$
\inf_{x \in A_n(a)} \sigma_x^{2d} \geq e^{-5\beta a}.
$$
Since $a$ may be as small as we wish, we conclude that necessarily
\be
\forall \ep>0 \ \exists\, \delta>0 \ \big[
A(x) \leq \delta \, \Rightarrow \, {\rm det} K''(x) \geq 1-\ep\big].
\en
But this is the same as (8.3), which is equivalent to (8.2).

Let us now explain why
\be
\lim_{A(x) \rightarrow 0} K''(x) = I_d \ \Longleftrightarrow \, 
\lim_{A(x) \rightarrow 0} {\rm det}\, K''(x) = 1.
\en
The implication $``\Rightarrow"$ is obvious.
For the opposite direction, we may assume that the necessary condition
(8.5) is fulfilled. Since ${\rm det} K''(x) = \lambda_1(x) \dots \lambda_d(x)$
in terms of the eigenvalues $\lambda_j(x)$ of $R_x$, this condition
may be stated in a weaker form as
\be
\forall \ep>0 \ \exists\, \delta>0 \ \big[
A(x) \leq \delta \, \Rightarrow \, \min_{1 \leq j \leq d} \ \lambda_j(x) \geq (1-\ep)^d\big].
\en
In order to reverse the conclusion, let us return to Corollary 4.2 and recall that 
the eigenvalues satisfy
$$
\limsup_{A(x) \rightarrow 0} \, \lambda_j(x) \leq 1, \quad 1 \leq j \leq d,
$$
for any $x \in \R^d$, which is a stronger property compared to (8.1). 
Thus, for any $\ep > 0$, there exists $\delta>0$ such that
$$
A(x) \leq \delta \, \Rightarrow \, \max_{1 \leq j \leq d} \lambda_j(x) \leq 1+\ep.
$$
Being combined with (8.7) with $0 < \ep < 1$, this gives
$(1-\ep)^d \leq \lambda_j(x) \leq 1+\ep$ and hence
$|\lambda_j(x) - 1| \leq 1 - (1-\ep)^d$ for all $j = 1,\dots,d$, as long as 
$A(x) \leq \delta$. In this case it follows that
$$
\|K''(x) - I_d\|_{\rm HS}^2 = \sum_{j=1}^d\, |\lambda_j(x) - 1|^2 \leq 
d\,\big(1 - (1-\ep)^d\big).
$$
As a result, we obtain the opposite implication  in (8.6).

It remains to see that the property 
\be
\lim_{A(x) \rightarrow 0} A''(x) = 0
\en
may be restated as the conditions $a)-b)$ in Theorem 1.2: 

\vskip3mm
$a)$ $A''(x) = 0$ for every point $x \in \R^d$ such that $A(x) = 0$;

$b)$ $\lim_{k \rightarrow \infty} A''(x_k) = 0$ for every sequence
$|x_k| \rightarrow \infty$ such that $A(x_k) \rightarrow 0$ as $k \rightarrow \infty$.

\vskip3mm
Obviously, the conditions $a)-b)$ follow from (8.8).
For the converse direction, assume that (8.8) is not true. Then
there would exist $\ep > 0$ such that, for any $\delta > 0$, one can pick up
a point $x \in \R^d$ with the property that
$$
A(x) \leq \delta \quad {\rm and} \quad \|A''(x)\|_{\rm HS} > \ep.
$$
Choosing $\delta = \delta_k \downarrow 0$, we would obtain a sequence $x_k \in \R^d$
such that $A(x_k) \leq \delta_k$ and $\|A''(x_k)\|_{\rm HS} > \ep$. If this sequence
is bounded, it would contain a convergent subsequence
$x_{k_l} \rightarrow x$ with $A(x) = 0$ and $\|A''(x)\|_{\rm HS} \geq \ep$,
by continuity of the  functions $A$ and $A''$. But this contradicts $a)$.
In the other case, one can subtract a subsequence such that
$|x_{k_l}| \rightarrow \infty$,  $A(x_{k_l}) \rightarrow 0$ as $l \rightarrow \infty$,
while $\|A''(x_{k_l})\|_{\rm HS} > \ep$. But this contradicts $b)$.
\qed

\vskip7mm
\section{{\bf Proof of Corollaries 1.3-1.4}}
\setcounter{equation}{0}

\vskip2mm
\noindent
As in Theorem 1.2, suppose that the random vector $X$ has mean 
zero and identity covariance matrix. In addition, assume that:

\vskip4mm
1) $Z_n$ has density $p_n$ for some $n=n_0$ such that 
$T_\infty(p_n||\varphi) < \infty$;

2) $X$ is strictly sub-Gaussian: $L(t) \leq e^{|t|^2/2}$
or equivalently $\Psi(t) \leq 1$ for all $t \in \R^d$.

\vskip4mm
\noindent
Recall that, by the separation property, we mean the relation
\be
L(t) \leq (1 - \delta)\, e^{|t|^2/2},
\en
which holds for all $t_0 > 0$ and $|t| \geq t_0$ with some 
$\delta = \delta(t_0)$, $\delta \in (0,1)$.

\vskip5mm
{\bf Proof of Corollary 1.3.}
From (9.1) it follows that the log-Laplace transform and the function $A$ satisfy
$$
K(t) \leq \frac{1}{2}\, |t|^2 - \log(1 - \delta), \quad A(t) \geq - \log(1 - \delta).
$$
Hence, the approach $A(t) \rightarrow 0$ is only possible when $t \rightarrow 0$.
But, for strictly sub-Gaus\-sian distributions, we necessarily have 
$A(t) = O(|t|^4)$ and $A''(t) = O(|t|^2)$ near zero.
Therefore, the condition (8.8) is fulfilled automatically.
\qed

\vskip4mm
Next, let us apply Theorem 1.2 to the Laplace transforms 
$L(t) = \E\,e^{\left<t,X\right>}$ with 
$$
\Psi(t) = L(t)\,e^{-|t|^2/2}, \quad t \in \R^d,
$$
being periodic, that is, satisfying 
\be
\Psi(t+h) = \Psi(t), \quad t \in \R^d.
\en
for some $h = (h_1,\dots,h_d) \in \R^d$, $h \neq 0$. Without loss of generality, let 
$h_i \geq 0$ (and not all of them are zero). Put 
$[0,h] = [0,h_1] \times \dots \times [0,h_d]$.

\vskip5mm
{\bf Proof of Corollary 1.4.} The function $\Psi(t)$ is positive, and
by the assumption (9.2), the function
$$
A(t) = -\log \Psi(t) = \frac{1}{2}\,|t|^2 - K''(t)
$$
is $h$-periodic as well. To express the condition $a)$ in terms of $\Psi$, 
let us differentiate the equality $\Psi(t) = e^{-A(t)}$ to get
$$
\partial_{t_i} \Psi(t) = -\partial_{t_i} A(t)\,e^{-A(t)}
$$
and
\bee
\partial^2_{t_i t_j} \Psi(t) 
 & = &
-\partial^2_{t_i t_j} A(t)\,e^{-A(t)} +
\partial_{t_i} A(t)\,\partial_{t_j} A(t)\,e^{-A(t)}.
\ene
Thus, $A''(t) = -\Psi''(t)$ for every point $t \in \R^d$ such $A(t) = 0$.
Recall that in this case, necessarily $A'(t) = 0$ and therefore $\Psi'(t) = 0$.

As for the condition $b)$ in Theorem 1.2, it may be reduced to $a)$. Indeed, 
assume that $A(x_k) \rightarrow 0$ for some sequence $x_k \in \R^d$
such that $|x_k| \rightarrow \infty$. Using the condition $a)$, we need to show 
that $A''(x_k) \rightarrow 0$. The latter is equivalent to the assertion that 
from any subsequence $x_{k_l}$ one may subtract a further subsequence $x_{k_l'}$
such that $A''(x_{k_l'}) \rightarrow 0$. For simpicity, let a given subsequence
be the whole sequence $x_k$. By the periodicity (9.2), $A(x_k) = A(y_k)$ and 
$A''(x_k) = A''(y_k)$ for some $y_k \in [0,h]$. By compactness, there is 
a convergence subsequence $y_{k_l} \rightarrow y \in [0,h]$ as $l \rightarrow \infty$.
But then, by continuity, $A(y) = 0$ and hence $A''(y)= 0$, by $a)$.
As a consequence, $A''(x_{k_l}) = A''(y_{k_l}) \rightarrow A''(y)= 0$.
\qed

\vskip7mm
\section{{\bf Laplace Transforms with Separation Property}}
\setcounter{equation}{0}

\vskip2mm
\noindent
In dimension $d=1$, Corollary 1.3 can be illustrated by different examples.
Let us recall two results from \cite{B-C-G4}, assuming that
$X$ is a sub-Gaussian random variable with mean zero. In this case 
the characteristic function
$$
f(z) = \E\,e^{izX}, \quad z \in \C,
$$
represents an entire function in the complex plane of order at most 2.
If $f(z)$ does not have any real or complex zeros, a well-known theorem 
due to Marcinkiewicz \cite{M} implies that the distribution of $X$ is already Gaussian.
Thus, nonnormal sub-Gaussion distributions have characteristic functions that
need to have zeros. 

The strict sub-Gaussianity is defined by the relation
\be
\E\,e^{tX} \leq e^{\sigma^2 t^2/2}, \quad t \in \R,
\en
where $\sigma^2 = \Var(X)$ is the variance of $X$.

\vskip5mm
{\bf Proposition 10.1.} {\sl If the distribution of $X$ is symmetric, and
all zeros of $f(z)$ with ${\rm Re}(z) \geq 0$ lie in the
cone centered on the real axis defined by
$$
|{\rm Arg}(z)| \leq \frac{\pi}{8},
$$ 
then $X$ is strictly sub-Gaussian. Moreover, if $X$ is nonnormal, then
for any $t_0 > 0$, there exists $c = c(t_0)$ in the interval $0 < c < \sigma^2$ 
such that
\be
\E\,e^{tX} \leq e^{c t^2/2}, \quad |t| \geq t_0. 
\en
}

The inequality (10.2) strengthens not only (10.1), but also the separation 
relation (9.1) (for $\sigma^2 = 1$). The claim about the strict sub-Gaussianity
in Proposition 10.1 refines a theorem due to Newman \cite{N}, who
considered the case where $f(z)$ has only real zeros (cf. also \cite{Bu-K}).
It was also shown in \cite{B-C-G4} that the condition $|{\rm Arg}(z)| \leq \frac{\pi}{8}$
is also necessary for the strict sub-Gaussianity of $X$, when it has
a symmetric distribution, and its characteristic function $f(z)$ has exactly one
zero $z$ in the quandrant ${\rm Re}(z) > 0$, ${\rm Im}(z) > 0$.

\vskip3mm
{\bf Proposition 10.2.} {\sl  If $X$ is nonnormal, and
the function $K(\sqrt{|t|})$ is concave on the half-axis 
$t>0$ and is concave on the half-axis $t<0$, then $(10.2)$ holds true.
}

\vskip5mm
In \cite{B-C-G4} one can find various examples illustrating these propositions.
In particular, the symmetric Bernoulli and the uniform distribution on 
a symmetric interval are strictly subgaussian, as well as convergent infinite 
convolutions of such distributions. Moreover, they satisfy the separation property (10.2).

Turning to the multidimensional situation, let us only mention two examples.
For a subgaussian random vector $X$ in $\R^d$ with mean zero and covariance
matrix $\sigma^2 I_d$, the notion of the strict subgaussianity is defined by
\be
\E\,e^{\left<t,X\right>} \leq e^{\sigma^2 |t|^2/2}, \quad t \in \R^d.
\en
This class of probability distributions is invariant under convolutions 
and weak limits. 

It should be clear that, a product measure on $\R^d$ satisfies (10.3), if and only 
if all marginals are strictly sub-Gaussian. For spherically invariant distributions, 
(10.3) is also reduced to dimension one.

\vskip3mm
{\bf Proposition 10.3.} {\sl Suppose that the distribution of a sub-Gaussian 
random vector $X = (X_1,\dots,X_d)$ is spherically invariant.
Then $X$ is strictly sub-Gaussian, if and only if $X_1$ is strictly sub-Gaussian.
In this case, if $X$ is nonnormal, and $X_1$ satisfies $(10.2)$, then
for any $t_0 > 0$, there exists $c = c(t_0)$,
$0 < c < \sigma^2$, such that
\be
\E\,e^{\left<t,X\right>} \leq e^{c |t|^2/2}, \quad t \in \R^d, \ |t| \geq t_0. 
\en
}

{\bf Proof.} By the assumption, the random vectors $X$ and $UX$ are 
equidistributed for any linear orthogonal transformation of the space $\R^d$. 
Given $t \in \R^d$, choose $U$
such that $U' t = |t| e_1 = |t|\,(1,0,\dots,0)$. Then
$$
\E\,e^{\left<t,X\right>} = \E\,e^{\left<t,UX\right>} = \E\,e^{\left<U't,X\right>} =
\E\,e^{|t| X_1}.
$$
Since the distribution of $X_1$ is symmetric about the origin, (10.3)
is equivalent to (10.1) with $X_1$ in place of $X$. The same is true
about the equivalence of (10.2) and (10.4).
\qed

\vskip5mm
Here are two basic examples, where as before, $p_n$ denotes the density of the
normalized sum of $n$ independent copies of $X$.

\vskip5mm
{\bf Corollary 10.4.} {\sl The uniform distribution on the Euclidean ball
$B(r)$ in $\R^d$ with center at the origin and radius $r>0$ satisfies $(10.3)$.
As a consequence, if  $r^2 = d+2$, then 
$T_\infty(p_n||\varphi) \rightarrow 0$ as $n \rightarrow \infty$.
}

\vskip5mm
The assumption $r^2 = d+2$ corresponds to the requirement that
$\E\,|X|^2 = d$ or $\E X_1^2 = 1$ for the random vector $X = (X_1,\dots,X_d)$
uniformly distributed in $B(r)$. This is equivalent to the property that $X$ 
has identity covariance matrix.

\vskip5mm
{\bf Corollary 10.5.} {\sl The same assertion holds for the uniform distribution 
on the Euclidean sphere in $\R^d$ with center at the origin and radius $r = \sqrt{d}$.
}

\vskip5mm
Note that the uniform distribution on the Euclidean sphere is not absolutely
continuous. But its $n$-th convolution power for large $n$ has a bounded,
compactly supported density, so that $T_\infty(p_n||\varphi) < \infty$
for some $n=n_0$.

Introduce one dimensional probability densities
$$
q_\alpha(x) = \frac{1}{c_\alpha}\, (1 - x^2)^{\alpha - 1}, \quad |x| < 1,
$$
with parameter $\alpha>0$, where (using the usual gamma-function)
\begin{eqnarray}
c_\alpha 
 & = &
\int_{-1}^1 (1 - x^2)^{\alpha - 1}\,dx \nonumber \\
 & = & 
\int_0^1 y^{-1/2}\,(1-y)^{\alpha-1}\,dy \, = \,
\frac{\Gamma(\frac{1}{2}) \Gamma(\alpha)}{\Gamma(\alpha + \frac{1}{2})}
\end{eqnarray}
is the normalizing constant. For the proof of Corollaries 10.4-10.5, we need:

\vskip5mm
{\bf Lemma 10.6.} {\sl For any $\alpha>0$, the distribution $\mu_\alpha$ 
with density $q_\alpha$ is strictly sub-Gaussian. Moreover, the separation 
property $(10.2)$ holds true.
}

\vskip5mm
{\bf Proof.} Expanding the cosh-function in power series, for the random variable 
$\xi_\alpha$ with density $q_\alpha$ the Laplace transform is given by
\bee
L_\alpha(t)
 & = &
\frac{1}{c_\alpha}\int_{-1}^1 e^{tx}\,(1 - x^2)^{\alpha - 1}\,dx \, = \,
\frac{2}{c_\alpha}\int_0^1 \cosh(tx)\,(1 - x^2)^{\alpha - 1}\,dx \\
 & = &
\frac{2}{c_\alpha}\, \sum_{n=0}^\infty\, \frac{t^{2n}}{(2n)!} \,
\int_0^1  x^{2n}\,(1 - x^2)^{\alpha - 1}\,dx \\
 & = &
\frac{1}{c_\alpha}\, \sum_{n=0}^\infty\, \frac{t^{2n}}{(2n)!} \,
\int_0^1  y^{n - \frac{1}{2}}\,(1 - y)^{\alpha - 1}\,dx \\
 & = &
\frac{1}{\sqrt{\pi}}\, \sum_{n=0}^\infty\, 
\frac{\Gamma(s)}{\Gamma(n+s)}\,
\frac{\Gamma(n + \frac{1}{2})}{(2n)!} \, t^{2n},
\ene
where we used (10.5) in the last step together with notation $s = \alpha + \frac{1}{2}$.
From this expansion we find that
$$
\sigma_\alpha^2 = \Var(\xi_\alpha) = \frac{1}{2\alpha + 1} = 
\frac{1}{2s} \quad (\sigma_\alpha>0).
$$
Hence
\bee
L_\alpha(t/\sigma_\alpha)
 & = &
\frac{1}{\sqrt{\pi}}\, \sum_{n=0}^\infty\, 
\frac{\Gamma(s) s^n}{\Gamma(n+s)}\,
\frac{\Gamma(n + \frac{1}{2})}{(2n)!} \,(2t^2)^n \\
 & = &
1 + \frac{1}{2}\,t^2 + \frac{1}{\sqrt{\pi}}\, \sum_{n=2}^\infty\, 
\frac{s^{n-1}}{(s+1) \dots (s + n-1)}\,
\frac{\Gamma(n + \frac{1}{2})}{(2n)!} \,(2t^2)^n.
\ene
The fraction inside the sum represents an increasing function in $s$, and
letting $s \rightarrow \infty$, we get
$$
L_\alpha(t/\sigma_\alpha) \leq \frac{1}{\sqrt{\pi}}\, \sum_{n=0}^\infty\, 
\frac{\Gamma(n + \frac{1}{2})}{(2n)!} \, (2t^2)^n = \cosh(t).
$$
Thus, the Laplace transform of $\xi_\alpha/\sigma_\alpha$ is bounded
by the Laplace transform of the Bernoulli distribution 
$\mu = \frac{1}{2}\,\delta_{-1} + \frac{1}{2}\,\delta_1$
(which is consistent with the property that $\mu_\alpha \rightarrow \mu$
weakly as $\alpha \rightarrow \infty$).
\qed

\vskip5mm
{\bf Proof of Corollaries 10.4-10.5.} 
Without loss of generality, let $r=1$. If the random vector
$X = (X_1,\dots,X_d)$ is uniformly distributed in the unit ball $B(1)$, 
the distribution of $X_1$ has density $q_\alpha$ with $\alpha = (d+1)/2$. 
Similarly, if $X$ is uniformly distributed in the unit sphere, the distribution 
of $X_1$ has density $q_\alpha$ with $\alpha = (d-1)/2$. It remains 
to apply Lemma 10.6 and Proposition 10.3.
\qed

\vskip7mm
\section{{\bf Laplace Transforms with Periodic Components }}
\setcounter{equation}{0}

\vskip2mm
\noindent
In order to describe examples illustrating Corollary 1.4, let us start with the following
definition. We write $h = (h_1,\dots,h_d) \geq 0$, if all $h_j \geq 0$.

\vskip5mm
{\bf Definition.} We say that the distribution $\mu$ of a random 
vector $X$ in $\R^d$ is periodic with respect to the standard normal law, 
with period $h \geq 0$ ($h \neq 0$), if it has a density $p(x)$ such that the function
$$
q(x) = \frac{p(x)}{\varphi(x)} = \frac{d\mu(x)}{d\gamma(x)},
\quad x \in \R^d,
$$
is periodic with period $h$, that is, $q(x+h) = q(x)$ for all $x \in \R^d$.

\vskip5mm
Here, $q$ represents the density of $\mu$ with respect to the standard Gaussian 
measure $\gamma$ on $\R^d$. We denote the class of all such distributions by 
$\mathfrak F_h$, and say that $X$ belongs to $\mathfrak F_h$.
In dimension $d=1$, this class was studied in \cite{B-C-G4}, and here we extend
a number of one dimensional observations to higher dimensions.
For this aim, we introduce the componentwise multiplication of vectors
$$
xy = (x_1 y_1,\dots,x_d y_d), \quad 
x = (x_1,\dots,x_d), \ y = (y_1,\dots,y_d) \in \R^d.
$$

\vskip4mm
{\bf Proposition 11.1.}  {\sl If $X$ belongs to the class $\mathfrak F_h$, then 
for all $m \in \Z^d$,
\be
\E\,e^{\left<mh,X\right>} = e^{|mh|^2/2}.
\en
In particular, the random vector $X$ is sub-Gaussian. 
}

\vskip5mm
{\bf Proof.}
By the periodicity, $q(x - mh) = q(x)$ for all $x \in \R^d$
and $m \in \Z^d$. Hence, the random vector $X + mh$ 
has density
\bee
p(x - mh) 
 & = &
q(x - mh) \varphi(x - mh) \\
 & = &
q(x)\varphi(x)\,e^{\left<mh,x\right> - \frac{1}{2}\,|m h|^2}
 \, = \,
p(x)\,e^{\left<mh,x\right> - \frac{1}{2}\,|m h|^2}.
\ene
It remains to integrate this equality over the variable $x$, which leads to (11.1).

Next, starting from (11.1), it is easy to see that $\E\,e^{c |X|^2} < \infty$
for some $c>0$.
\qed

 \vskip5mm
As a consequence, the Laplace transform $L(t) = \E\,e^{\left<t,X\right>}$, 
$t \in \R^d$, is finite and may be extended to the $d$-dimensional complex 
space $\C^d$ as an entire function $L(z) = L(z_1,\dots,z_d)$. 
By saying
``entire", it is meant that a given function on $\C^d$ is entire with respect 
to every complex coordinate $z_j = t_j + i y_j$ ($t_j,y_j \in \R$)
and is $C^\infty$-smooth as a functon of $2d$ real variables
$t_1,y_1,\dots,t_d,y_d$. 

This property of $L(z)$ may be further refined.

\vskip5mm
{\bf Proposition 11.2.}  {\sl If $X$ belongs to $\mathfrak F_h$, then
its Laplace transform is an entire function of order 2 (with respect to every
complex coordinate). Moreover,
\be
|L(z)| \leq e^{|t|^2 + |h|^2}, \quad z = t+iy \in \C^d.
\en
}

{\bf Proof.} 
Since ${\rm Re} \left<z,X\right> = \left<t,X\right>$ for $z = t+iy$,
we have $|L(z)| \leq L(t)$. Hence, one may assume that $y = 0$ in (11.2).

We employ the convexity of the function $K(t) = \log L(t)$.
For simplicity, suppose that $t = (t_1,\dots,t_d) \in \R^d_+$.
Take an integral vector $m = (m_1,\dots,m_d)$ with positive components 
such that $(m_j - 1) h_j \leq t_j < m_j h_j$ for all $j \leq d$. Since
$t$ lies in the cube with sides $[0,m_j h_j]$, this point may be written as
a convex mixture of vertices of the cube
$$
t = \sum_\ep a_\ep v_\ep, \quad a_\ep \geq 0, \ \ \sum_\ep a_\ep = 1.
$$
Here the vertex $v_\ep = (\ep_1 m_1 h_1,\dots,\ep_d m_d h_d)$ is 
parametrized by the tuple $\ep = (\ep_1,\dots,\ep_d)$ with $\ep_j = 0$ 
or $\ep_j = 1$ for each $j \leq d$. By Jensen's inequality and (11.1),
\bee
K(t) 
 & \leq &
\sum_\ep a_\ep K(v_\ep) \, = \, \frac{1}{2} \sum_\ep a_\ep
\bigg[\sum_{j=1}^d \ep_j^2 m_j^2 h_j^2\bigg] \\
 & \leq & 
\frac{1}{2} \sum_{j=1}^d m_j^2 h_j^2 \, \leq \, 
\frac{1}{2} \sum_{j=1}^d\, (t_j + h_j)^2.
\ene
Dropping the condition on the sign of $t_j$, more generally we obtain that
$$
K(t) \leq \frac{1}{2} \sum_{j=1}^d\, (|t_j| + h_j)^2 \leq 
\sum_{j=1}^d\, (t_j^2 + h_j^2) = |t|^2 + |h|^2.
$$

Thus, we obtain (11.2) which shows that $L(z)$ is an entire function 
of order at most 2. On the other hand, by (11.1), it is
is an entire function of order at least 2.
\qed

\vskip5mm
Let us also mention the periodicity property for convolutions.

\vskip5mm
{\bf Proposition 11.3.} {\sl If $X$ belongs to $\mathfrak F_h$, then
$Z_n$ belongs to $\mathfrak F_{h\sqrt{n}}$.
}

\vskip5mm
Here as before, $Z_n = \frac{1}{\sqrt{n}}\,(X_1 + \dots + X_n)$ denotes 
the normalized sum of independent copies $X_k$ of the random vector $X$.
The proof is similar to the proof of Proposition 10.4 from \cite{B-C-G4}
for the one dimensional case, so we omit it.

\vskip7mm
\section{{\bf Characterizations of Periodicity in terms of Laplace Transform}}
\setcounter{equation}{0}

\vskip2mm
\noindent
Fix $h = (h_1,\dots,h_d) \in \R^d$, $h \neq 0$, for simplicity with $h_j \geq 0$.
Here we prove:

\vskip5mm
{\bf Proposition 12.1.} {\sl If $X$ belongs to $\mathfrak F_h$, then
the function
\be
\Psi(t) = L(t)\,e^{-|t|^2/2}, \quad t \in \R^d, 
\en
is periodic with period $h$. It can be extended to the complex space 
$\C^d$ as an entire function of order at most $2$. Conversely, if the
fuction $\Psi(t)$ for a sub-Gaussian random vector $X$ is $h$-periodic, then 
$X$ belongs to $\mathfrak F_h$, as long as 
the characteristic function $f(t)$ of $X$ is integrable on $\R^d$. 
}

\vskip5mm
The last claim is based on the following general observation.

\vskip5mm
{\bf Lemma 12.2.} {\sl Let  $f(z) = \E\,e^{i\left<z,X\right>}$ $(z \in \C^d)$
be the characteristic function of a sub-Gaussian random vector $X$. 
If $f(t)$ and $f(t+ih)$ are integrable in $t \in \R^d$, then
\be
\int_{\R^d} e^{-i \left<t,x\right>} f(t)\,dt = 
\int_{\R^d} e^{-i \left<t + ih,x\right>} f(t + ih)\,dt.
\en
}

\vskip2mm
{\bf Proof.} Write $z = (z_1,\dots,z_d)$ with $z_j \in \C$.
By assumption, $f(z) = f(z_1,\dots,z_d)$ is well-defined and represents
an entire function. In addition,
\be
|f(z)| \leq \E\,|e^{i\left<z,X\right>}| = \E\,e^{-\left<y,X\right>} \leq
\E\,e^{|y|\,|X|}, \quad z = t+iy, \ t,y \in \R^d.
\en

One may rewrite the first integral in (12.2) in a different way using 
contour integration. For this aim, let $X_\ep = X + \ep Z$, where 
$\ep>0$ and $Z$ is a standard normal random vector in $\R^d$ 
independent of $X$. The random vector $X_\ep$ is also sub-Gaussian 
and has an entire characteristic function
$$
f_\ep(z) = f(z)\,e^{- \ep^2 z^2/2}, \quad z \in \C^d,
$$
where $z^2 = z_1^2 + \dots + z_d^2$.
Moreover, by (12.3), for all $t \in \R^d$,
\be
|f_\ep(t + iy)| \leq K e^{-\ep^2 |t|^2/2}, \quad y = (y_1,\dots,y_j), \ |y_j| \leq h_j,
\en
with some constant $K$ which does not depend on $t$.

Putting $x = (x_1,\dots,x_d)$, $t = (t_1,\dots,t_d)$, let us integrate 
the left integrand in (12.2) with respect to $t_1$, keeping the remaining
variables $t_2,\dots,t_d$ fixed.
Given $T>0$, consider the rectangle contour with sides
\bee
C_1 & = &
[-T,T], \qquad \qquad \qquad C_2 = [T,T+ih_1], \\ 
C_3 & = &
[T+ih_1,-T+ih_1], \quad C_4 = [-T+ih_1,-T],
\ene
so that to apply Cauchy's theorem which yields
\bee
\int_{C_1} e^{-iz_1 x_1} f_\ep(z)\,dz_1 + \int_{C_2} e^{-iz_1 x_1} f(z)\,dz_1 
 & & \\ 
 & & \hskip-40mm + \ 
\int_{C_3} e^{-iz_1 x_1} f_\ep(z)\,dz_1 + \int_{C_4} e^{-iz_1 x_1} f_\ep(z)\,dz_1 = 0.
\ene
For points $z_1 = t_1 + iy_1$ on the contour, 
$|e^{-iz_1 x_1}| = e^{x_1 y_1} \leq e^{|x| h_1}$. Hence, by the decay property
(12.4), the integrals over $C_2$ and $C_4$ are vanishing as $T \rightarrow \infty$,
and
\bee
\int_{-\infty}^\infty e^{-it_1 x_1} f_\ep(t)\,dt_1
 & = &
\lim_{T \rightarrow \infty}\, \int_{C_1} e^{-iz_1 x_1} f_\ep(z)\,dz_1 \\
 & = & - \ 
\lim_{T \rightarrow \infty}\, \int_{C_3} e^{-iz_1 x_1} f_\ep(z)\,dz_1 \\
 & = &
\int_{-\infty}^\infty e^{-i(t_1 + ih_1) x_1} f_\ep(t_1+ih_1,t_2,\dots,t_d)\,dt_1.
\ene
Thus,
\be
\int_{-\infty}^\infty e^{-it_1 x_1} f_\ep(t)\,dt_1 = 
\int_{-\infty}^\infty e^{-i(t_1 + ih_1) x_1} f_\ep(t_1+ih_1,t_2,\dots,t_d)\,dt_1.
\en

Turning to the next variable $t_2$, first note that the function
$$
g(z_2,\dots,z_d) = \int_{-\infty}^\infty e^{-i(t_1 + ih_1) x_1} 
f_\ep(t_1+ih_1,z_2,z_3,\dots,z_d)\,dt_1, \quad z_j = t_j + iy_j \in \C,
$$
is entire and admits a similar bound as (12.4)
$$
|g(z_2,\dots,z_d)| \leq K\,e^{-\ep^2 (t_2^2 + \dots + t_d^2)/2}, \quad |y_j| \leq h_j.
$$
Hence, one may perform the contour integration like in the previois step
leading to
$$
\int_{-\infty}^\infty e^{-it_2 x_2}\, g(t_2,t_3,\dots,t_d)\,dt_2 = 
\int_{-\infty}^\infty e^{-i(t_2 + ih_2) x_2}\, g(t_2+ih_2,t_3,\dots,t_d)\,dt_2.
$$
By the Fubini theorem and using (12.5), we then get
\bee
\int_{-\infty}^\infty \int_{-\infty}^\infty 
e^{-it_1 x_1 - it_2 x_2} f_\ep(t)\,dt_1 dt_2 
 & & \\
 & & \hskip-50mm = \ 
\int_{-\infty}^\infty \int_{-\infty}^\infty 
e^{-i(t_1 + ih_1) x_1 - i(t_2 + ih_2) x_2} 
f_\ep(t_1+ih_1,t_2+ih_2,t_3,\dots,t_d)\,dt_1 dt_2.
\ene

Repeating the process, we will arive at the equality for the $d$-dimensional integrals
$$
\int_{\R^d} e^{-i \left<t,x\right>} f_\ep(t)\,dt = 
\int_{\R^d} e^{-i \left<t + ih,x\right>} f_\ep(t + ih)\,dt.
$$
It remains to let $\ep \rightarrow 0$ and make use of the Lebesgue
dominated convergence theorem, which provides the desired equality in (12.2).
\qed

\vskip5mm
{\bf Proof of Proposition 12.1.} By periodicity of $q$, changing the variable 
$x = y+h$, we have, for any $t \in \R^d$,
\bee
L(t+h) 
 & = &
\int_{\R^d} e^{\left<t+h,x\right>}\,q(x)\,\varphi(x)\,dx \\
 & = &
\int_{\R^d} e^{\left<t+h,y+h\right>}\,q(y+h)\,\varphi(y+h)\,dy \\
 & = &
\int_{\R^d} e^{\left<t+h,y+h\right>}\,q(y)\,\varphi(y)\,
e^{-\left<y,h\right> - |h|^2/2}\,dy
 \, = \,
L(t)\,e^{\left<t,h\right> + |h|^2/2}.
\ene
Hence
$$
L(t+h) \,e^{-|t+h|^2/2} = L(t)\,e^{-|t|^2/2},
$$ 
which is the first claim of the proposition. Since $L(z)$ is an entire 
function of order 2,
the formula (12.1) admits a natural extension to $\C^d$
$$
\Psi(z) = L(z)\,e^{-\frac{1}{2}\,(z_1^2 + \dots +z_d^2)}, \quad 
z = (z_1,\dots,z_d) \in \C^d,
$$
which is an entire function with respect to every component $z_j$
of order at most 2. By analyticity and periodicity on $\R^d$,
\be
\Psi(z+h) = \Psi(z) \quad {\rm for \ all} \ z \in \C^d.
\en

It remains to prove the periodicity of the density $q(x) = p(x)/\varphi(x)$. 
The characteristic function 
of $X$ admits an entire extension using the formula
$$
f(z) = \E\,e^{i\left<z,X\right>} = L(iz) = \Psi(iz)\,e^{-z^2/2}, \quad z \in \C^d,
$$
where again $z^2 = z_1^2 + \dots + z_d^2$. Hence, by (12.6),
$$
f(t+ih)\,e^{(t+ih)^2/2} = f(t)\,e^{t^2/2},
$$
that is,
\be
f(t+ih) = f(t)\,e^{-i \left<t,h\right> + |h|^2/2} \quad 
{\rm for \ all} \ t \in \R^d.
\en
This identity also shows that, due to the integrability of $f(t)$, the function
$f(t+ih)$ is integrable as well.

Again by the the integrability of $f(t)$, the random vector $X$ 
has a continuous density given by the Fourier inversion formula
$$
p(x) = 
\frac{1}{(2\pi)^d} \int_{\R^d} e^{-i\left<t,x\right>} f(t)\,dt, \quad 
x \in \R^d.
$$
This yields
$$
q(x) = \frac{p(x)}{\varphi(x)} = 
\frac{1}{(2\pi)^{d/2}}\,e^{|x|^2/2} \int_{\R^d} e^{-i\left<t,x\right>} f(t)\,dt
$$
and
$$
q(x+h) = \frac{1}{(2\pi)^{d/2}}\,e^{|x|^2/2}\, e^{\left<x,h\right> + |h|^2/2} 
\int_{\R^d} e^{-i\left<t,x\right> - i \left<t,h\right>} f(t)\,dt.
$$
Hence, we need to show that
$$
\int_{\R^d} e^{-i \left<t,x\right>} f(t)\,dt = e^{\left<x,h\right> + |h|^2/2} 
\int_{\R^d} e^{-i \left<t,x\right> - i \left<t,h\right>} f(t)\,dt.
$$
Here, the left integral may be rewritten according to Lemma 12.1, so that
the above equality is restated as
\be
\int_{\R^d} e^{-i \left<t + ih,x\right>} f(t + ih)\,dt = e^{\left<x,h\right> + |h|^2/2} 
\int_{\R^d} e^{-i \left<t,x\right> - i \left<t,h\right>} f(t)\,dt.
\en
Moreover, by (12.7), the first integrand is equal to
$$
e^{-i\left<t + ih,x\right>}\, \,e^{-i\left<t,h\right> + |h|^2/2}\,f(t).
$$
This proves (12.8).
\qed

\vskip4mm
{\bf Remark.} Since $f(t) = L(it) = \Psi(it)\,e^{-t^2/2}$, $t \in \R^d$,
the integrability assumption in Proposition 12.1 is fulfilled,
if $\Psi(z)$ has order smaller than 2.

\vskip7mm
\section{{\bf Periodic Components via Trigonometric Series}}
\setcounter{equation}{0}

\vskip2mm
\noindent
Proposition 12.1 is applicable to a variety of interesting examples including
the underlying distributions whose Laplace transform has the form 
\be
L(t) = \Psi(t)\, e^{|t|^2/2}, \quad t \in \R^d,
\en
where $\Psi$ is a $2\pi$-periodic function of the form
\be
\Psi(t) = 1 - c P(t), \quad P(t) = \sum_{k \in \Z^d} c_k e^{i\left<k,t\right>}
\en
(for simplicity, in the sequel we say ``$2\pi$-periodic" instead of
``$(2\pi,\dots,2\pi)$-periodic").
Here $c_k = a_k - ib_k$ are complex coefficients which are supposed to satisfy
\be
\sum_{k \in \Z^d}^\infty e^{|k|^2/2}\,|c_k| < \infty,
\en
and $c \in \R$ is a non-zero parameter. To ensure that $P(t)$ is real-valued,
we assume that $a_{-k} = a_k$ and $b_{-k} = -b_k$ for all $k \in \Z^d$,
in which case
$$
P(t) = \sum_{k \in \Z^d} 
\big(a_k\cos\left<k,t\right> + b_k\sin\left<k,t\right>\big).
$$
Any such function with real coefficients $a_k$ and $b_k$ can be written 
in the form (13.2).

Let us note that the function $\Psi(t)$ defined in the equality (13.1) for 
a sub-Gaussian random vector $X$ is smooth with $\Psi(0) = L(0) = 1$.
Differentiating (13.1), we see that $X$ has mean zero and identity
covariance matrix if and only if $\Psi'(0) = 0$ in vector sense and $\Psi''(0) = 0$ 
in matrix sense. For the $\Psi$-functions as in (13.2), this is equivalent to 
$P(0) = P'(0) = P''(0) = 0$, that is,
$$
\sum_{k \in \Z^d} a_k = 0, \quad \sum_{k \in \Z^d} b_k k = 0, \quad
\sum_{k=1}^\infty a_k\, k \otimes k = 0,
$$
where $k \otimes k$ denotes the $d \times d$ matrix with entries $k_i k_j$,
$1 \leq i,j \leq d$, $k = (k_1,\dots,k_d)$.
All the series are well convergent due to the condition (13.3).

\vskip5mm
{\bf Proposition 13.1.} {\sl If $P(0) = P'(0) = P''(0) = 0$ and $|c|$ is small 
enough, then $L(t)$ represents the Laplace transform of a sub-Gaussian 
random vector $X$ with mean zero, identity covariance matrix, and with density 
$p = q \varphi$, where $q$ is a bounded, $2\pi$-periodic function. This 
random vector is strictly sub-Gaussian, if $P(t) \geq 0$ for all $t \in \R^d$ 
and if $c > 0$ is small~enough.
}

\vskip5mm
{\bf Proof.} The functions
$u_\lambda(x) = \varphi(x)\cos\left<\lambda,x\right>$ and 
$v_\lambda(x) = \varphi(x)\sin\left<\lambda,x\right>$ with 
parameter $\lambda \in \R^d$ have respectively the Laplace transforms
\bee
\int_{\R^d} e^{\left<t,x\right>} u_\lambda(x)\,dx
 & = &
e^{-|\lambda|^2/2}\,\cos\left<\lambda,t\right> e^{|t|^2/2}, \\
\int_{\R^d} e^{\left<t,x\right>} v_\lambda(x)\,dx
 & = &
e^{-|\lambda|^2/2}\,\sin\left<\lambda,t\right> e^{|t|^2/2}.
\ene
Define
\be
q(x) = 1 - c \sum_{k \in \Z^d} 
e^{|k|^2/2}\,\big(a_k\cos\left<k,x\right> + b_k\sin\left<k,x\right>\big).
\en
It is bounded due to the condition (13.3) and is non-negative for sufficiently
small $|c|$, more pricesely, if
$$
\sum_{k \in \Z^d}^\infty e^{|k|^2/2}\,(|a_k| + |b_k|) \leq \frac{1}{|c|}.
$$
Moreover, the Laplace transform of the function 
$p(x) = q(x) \varphi(x)$ is exactly
$$
\int_{\R^d} e^{\left<t,x\right>} p(x)\,dx = (1 - cP(t))\,e^{-|t|^2/2},
\quad t \in \R^d.
$$

Recall that the requirement $P(0) = 0$ guarantees that $\int_{\R^d} p(x)\,dx = 1$,
while the property that $X$ has mean zero and identity covariance matrix
is equivalent to $P'(0) = P''(0) = 0$.

Finally, if $P(t) \geq 0$ for all $t \in \R^d$ and $c > 0$ is small enough,
then $0 < \Psi(t) \leq 1$, which means that $X$ is strictly sub-Gaussian.
\qed

\vskip7mm
\section{{\bf Examples Involving Trigonometric Polynomials}}
\setcounter{equation}{0}

\vskip2mm
\noindent
Some specific examples in Proposition 13.1 are based on trigonometric
polynomials. More precisely, suppose that the random vector $X$
has the Laplace transform
\be
L(t) = \Psi(t)\, e^{|t|^2/2}, \quad t \in \R^d,
\en
with
\be
\Psi(t) = 1 - c Q(t)^2, \quad Q(t) = \sum_{k \in \Z^d} 
\big(a_k\cos\left<k,t\right> + b_k\sin\left<k,t\right>\big),
\en
where $c>0$, assuming that the sum contains only finitely many real
coefficients.

\vskip5mm
{\bf Corollary 14.1.} {\sl If $Q(0) = Q'(0) = 0$ and $c>0$ is small 
enough, $L(t)$ represents the Laplace transform of a strictly sub-Gaussian 
random vector $X$ with mean zero, identity covariance matrix, and with density 
$p = q \varphi$, where $q$ is a bounded, $2\pi$-pe\-r\-i\-o\-dic function.
}

\vskip5mm
Let us complement this statement with the assertion about the central limit
theorem for the distance
$$
T_\infty(p_n||\varphi) = \sup_x\, \frac{p_n(x) - \varphi(x)}{\varphi(x)},
$$
where $p_n$ denotes the density of the normalized sum $Z_n$ constructed
for $n$ independent copies of $X$.

\vskip3mm
{\bf Corollary 14.2.} {\sl Under the assumptions of Corollary $14.1$, 
$T_\infty(p_n||\varphi) \rightarrow 0$ as $n \rightarrow \infty$, if and only if
\be
\forall\, t \in [0,2\pi]^d \ \ \big[Q(t) = 0 \, \Rightarrow \, Q'(t) = 0\big].
\en
This is the case where $Q(t) \geq 0$ for all $t \in [0,2\pi]^d$.
}

\vskip5mm
{\bf Proof.} We apply Proposition 13.1 with $P(t) = Q(t)^2$. Then
$P(0) = Q(0)^2 = 0$ and $P'(0) = 2 Q(0)\, Q'(0) = 0$. In addition,
\be
P''(t) = 2\, Q(t) Q''(t) + 2\,Q'(t) \otimes Q'(t),
\en
so that $P''(0) = 0$ as well. This proves Corollary 14.1.

For the assertion in Corollary 14.2, we appeal to Corollary 1.4 and simplify the
implication $\Psi(t) = 0 \, \Rightarrow \, \Psi''(t) = 0$. Here, the
hypothesis means that $Q(t) = 0$, in which case
$\Psi''(t) = -c P''(t) = -2c\,Q'(t) \otimes Q'(t)$, according to (14.4).
The latter matrix is zero, if and only if $Q'(t) = 0$.
\qed

\vskip5mm
{\bf Example 14.3.} Consider the Laplace transform
$$
L(t) = \big(1 - c \sin^{2m}(t_1 + \dots + t_d)\big)\, e^{|t|^2/2}, \quad 
t = (t_1,\dots,t_d) \in \R^d,
$$
with a fixed integer $m \geq 2$, where $c>0$ is small enough 
(depending on $m$ and $d$).
In this case, the conditions of Corollary 14.1 are fulfilled for
$$
Q(t) = \sin^{m}(t_1 + \dots + t_d).
$$
Since the condition (14.3) in Corollary 14.2 is fulfilled as well, we obtain 
the assertion about the CLT (although $Q(t)$ does not need be non-negative
for odd values of $m$).

An interesting feature of this example is that $L(t) = e^{|t|^2/2}$ on
countably many hyperplanes $t_1 + \dots + t_d = \pi l$, $l \in \Z$.

\vskip5mm
{\bf Example 14.4.} Modifying the previous example, put
$$
P(t) = Q(t)^2, \quad Q(t) = 
\big(1 - 4\sin^2(t_1 + \dots + t_d)\big)  \sin^2(t_1 + \dots + t_d).
$$
In this case, $Q(0) = Q'(0) = 0$, so that the conditions of Corollary 14.1
are fulfilled.

We have $Q(t) = 0$ for $t_1 = \pi/6$, $t_j = 0$ for $j \geq 2$. At this point
$\partial_{t_1} Q(t) \neq 0$, so that the condition (14.3) is not fulfilled.
Hence, by Corollary 14.2, the CLT with respect to $T_\infty$ does not hold 
in this example.

\vskip5mm
{\bf Acknowledgement.}
The research has been supported by the NSF grant DMS-2154001 and 
the GRF -- SFB 1283/2 2021 -- 317210226.


\vskip5mm
\begin{thebibliography}{BH3}
\itemsep=4pt

\bibitem{A-B-B-N} 
Artstein, S.; Ball, K. M.; Barthe, F.; Naor, A. On the rate of convergence in 
the entropic central limit theorem. Probab. Theory Related Fields 129 (2004), 
no. 3, 381--390. 

\bibitem{Bar} 
Barron, A. R. Entropy and the central limit theorem. 
Ann. Probab. 14 (1986), no. 1, 336--342. 

\bibitem{B-C}
Bobkov, S. G.; Chistyakov, G. P. Bounds for the maximum of the density of the 
sum of independent random variables. (Russian) Zap. Nauchn. Sem. S.-Peterburg. 
Otdel. Mat. Inst. Steklov. (POMI) 408, Veroyatnost i Statistika. 18 (2012), 62--73, 324; 
transl. in J. Math. Sci. (N.Y.) 199 (2014), no. 2, 100--106.

\bibitem{B-C-G1}
Bobkov, S. G.; Chistyakov, G. P.; G\"otze, F.  Rate of convergence and 
Edgeworth-type expansion in the entropic central limit theorem. 
Ann. Probab. 41 (2013), no. 4, 2479--2512.

\bibitem{B-C-G2}
Bobkov, S. G.; Chistyakov, G. P.; G\"otze, F. Berry-Esseen bounds in 
the entropic central limit theorem. Probab. Theory Related Fields 
159 (2014), no. 3-4, 435--478.

\bibitem{B-C-G3}
Bobkov, S. G.; Chistyakov, G. P.; G\"otze, F. R\'enyi divergence and the central limit 
theorem. Ann. Probab. 47 (2019), no. 1, 270--323. 

\bibitem{B-C-G4}
Bobkov, S. G.; Chistyakov, G. P.; G\"otze, F. Strictly subgaussian distributions.
Electron. J. Probab. 29 (2024), article no. 62, 1--28.

\bibitem{B-G1}
Bobkov, S. G.; G\"otze, F. Central limit theorem for R\'enyi divergence of
infinite order. Preprint (2024). To appear in: Ann. Probab.

\bibitem{B-G2}
Bobkov, S. G.; G\"otze, F. Berry-Esseen bounds in the local limt theorem.
Preprint (2024).

\bibitem{B-M}
Bobkov, S. G.; Madiman, M. The entropy per coordinate of a random vector is 
highly constrained under convexity conditions. IEEE Trans. Inform. Theory 
57 (2011), no. 8, 4940--4954.

\bibitem{Bu-K}
Buldygin, V. V.; Kozachenko, Yu. V. Metric characterization of random variables 
and random processes. Translated from the 1998 Russian original by V. Zaiats. 
Transl. Math. Monogr., 188
American Mathematical Society, Providence, RI, 2000. xii+257 pp.


\bibitem{D} 
Daniels, H. E. Saddlepoint approximations in statistics. 
Ann. Math. Statist. 25 (1954), 631--650.

\bibitem{E}
Esscher, F. On the probability function in the collective theory of risk. 
Skandinavisk Aktuarietidskrift. 15 (3) (1932), 175--195.

\bibitem{K}
Khinchin, A. I. Mathematical Foundations of Statistical Mechanics.
Translated by G. Gamow. Dover Publications, Inc., New York, N.Y., 1949. viii+179 pp.

\bibitem{L}
Linnik, J. V. An information-theoretic proof of the central limit theorem with 
Lindeberg conditions. Theory Probab. Appl. 4 (1959), 288--299.

\bibitem{M}
Marcinkiewicz, J. J. Sur une propri\'et\'e de la loi de Gauss. (French) 
           Math. Z. 44 (1939), no. 1, 612--618.

\bibitem{N}
Newman, C. M. Inequalities for Ising models and field theories which obey the Lee-Yang theorem.
Comm. Math. Phys. 41 (1975), 1--9.

\bibitem{Pr}
Prokhorov, Yu. V. A local theorem for densities. (Russian) 
Doklady Akad. Nauk SSSR (N.S.) 83 (1952), 797--800. 


\end{thebibliography}
\end{document}